\documentclass[onecolumn]{IEEEtran}
\usepackage[latin1]{inputenc}
\usepackage[english]{babel}
\usepackage[margin=1in]{geometry}

\usepackage{blindtext}
\usepackage{amssymb}
\usepackage{hyperref}
\hypersetup{
    colorlinks=true,
    linkcolor=blue,
    filecolor=magenta,      
    urlcolor=cyan,
}
 
\usepackage{multirow}
\usepackage{makecell}
\usepackage{amsthm}
\usepackage[normalem]{ulem}
\usepackage{amsmath}
\usepackage{booktabs}
\usepackage{array}
\usepackage{cite}
\usepackage{soul}
\usepackage{mathtools}
\usepackage{bm}
\newtheorem{theorem}{Theorem}[section]
\newtheorem{corollary}{Corollary}[section]
\newtheorem{lemma}[theorem]{Lemma}
\newtheorem{proposition}[theorem]{Proposition}

\newtheorem{remark}{Remark}[section]

\newcommand\nc\newcommand
\nc{\cA}{\mathcal{A}}\nc{\cB}{\mathcal{B}}\nc{\cC}{\mathcal{C}}\nc{\cD}{\mathcal{D}}
\nc{\cE}{\mathcal{E}}\nc{\cF}{\mathcal{F}}\nc{\cG}{\mathcal{G}}\nc{\cH}{\mathcal{H}}
\nc{\cI}{\mathcal{I}}\nc{\cJ}{\mathcal{J}}\nc{\cK}{\mathcal{K}}\nc{\cL}{\mathcal{L}}
\nc{\cM}{\mathcal{M}}\nc{\cN}{\mathcal{N}}\nc{\cO}{\mathcal{O}}\nc{\cP}{\mathcal{P}}
\nc{\cQ}{\mathcal{Q}}\nc{\cR}{\mathcal{R}}\nc{\cS}{\mathcal{S}}\nc{\cT}{\mathcal{T}}
\nc{\cU}{\mathcal{U}}\nc{\cV}{\mathcal{V}}\nc{\cW}{\mathcal{W}}\nc{\cX}{\mathcal{X}}
\nc{\cY}{\mathcal{Y}}\nc{\cZ}{\mathcal{Z}}

\nc{\bba}{\mathbf{a}}\nc{\bbb}{\mathbf{b}}\nc{\bbc}{\mathbf{c}}\nc{\bbd}{\mathbf{d}}
\nc{\bbe}{\mathbf{e}}\nc{\bbf}{\mathbf{f}}\nc{\bbg}{\mathbf{g}}\nc{\bbh}{\mathbf{h}}
\nc{\bbi}{\mathbf{i}}\nc{\bbj}{\mathbf{j}}\nc{\bbk}{\mathbf{k}}\nc{\bbl}{\mathbf{l}}
\nc{\bbm}{\mathbf{m}}\nc{\bbn}{\mathbf{n}}\nc{\bbo}{\mathbf{o}}\nc{\bbp}{\mathbf{p}}
\nc{\bbq}{\mathbf{q}}\nc{\bbr}{\mathbf{r}}\nc{\bbs}{\mathbf{s}}\nc{\bbt}{\mathbf{t}}
\nc{\bbu}{\mathbf{u}}\nc{\bbv}{\mathbf{v}}\nc{\bbw}{\mathbf{w}}\nc{\bbx}{\mathbf{x}}
\nc{\bby}{\mathbf{y}}\nc{\bbz}{\mathbf{z}}

\nc{\bbA}{\mathbf{A}}\nc{\bbB}{\mathbf{B}}\nc{\bbC}{\mathbf{C}}\nc{\bbD}{\mathbf{D}}
\nc{\bbE}{\mathbf{E}}\nc{\bbF}{\mathbf{F}}\nc{\bbG}{\mathbf{G}}\nc{\bbH}{\mathbf{H}}
\nc{\bbI}{\mathbf{I}}\nc{\bbJ}{\mathbf{J}}\nc{\bbK}{\mathbf{K}}\nc{\bbL}{\mathbf{L}}
\nc{\bbM}{\mathbf{M}}\nc{\bbN}{\mathbf{N}}\nc{\bbO}{\mathbf{O}}\nc{\bbP}{\mathbf{P}}
\nc{\bbQ}{\mathbf{Q}}\nc{\bbR}{\mathbf{R}}\nc{\bbS}{\mathbf{S}}\nc{\bbT}{\mathbf{T}}
\nc{\bbU}{\mathbf{U}}\nc{\bbV}{\mathbf{V}}\nc{\bbW}{\mathbf{W}}\nc{\bfX}{\mathbf{X}}
\nc{\bbY}{\mathbf{Y}}\nc{\bbZ}{\mathbf{Z}}

\nc{\sA}{\mathsf{A}}\nc{\sB}{\mathsf{B}}\nc{\sC}{\mathsf{C}}\nc{\sD}{\mathsf{D}}
\nc{\sE}{\mathsf{E}}\nc{\sF}{\mathsf{F}}\nc{\sG}{\mathsf{G}}\nc{\sH}{\mathsf{H}}
\nc{\sI}{\mathsf{I}}\nc{\sJ}{\mathsf{J}}\nc{\sK}{\mathsf{K}}\nc{\sL}{\mathsf{L}}
\nc{\sM}{\mathsf{M}}\nc{\sN}{\mathsf{N}}\nc{\sO}{\mathsf{O}}\nc{\sP}{\mathsf{P}}
\nc{\sQ}{\mathsf{Q}}\nc{\sR}{\mathsf{R}}\nc{\sS}{\mathsf{S}}\nc{\sT}{\mathsf{T}}
\nc{\sU}{\mathsf{U}}\nc{\sV}{\mathsf{V}}\nc{\sW}{\mathsf{W}}\nc{\sX}{\mathsf{X}}
\nc{\sY}{\mathsf{Y}}\nc{\sZ}{\mathsf{Z}}

\newcommand{\abs}[1]{\left|#1\right|}
\newcommand{\ceilenv}[1]{\left\lceil #1 \right\rceil}
\newcommand{\floorenv}[1]{\left\lfloor #1 \right\rfloor}
\newcommand{\parenv}[1]{\left( #1 \right)}

\nc{\set}[1]{\llbracket #1 \rrbracket}

\newcommand{\bal}[1]{\begin{align}\label{#1}}
\newcommand{\eal}{\end{align}}

\renewcommand{\le}{\leqslant}
\renewcommand{\leq}{\leqslant}

\renewcommand{\geq}{\geqslant}

\renewcommand{\Bbb}{\mathbb}

\newcommand{\eqdef}{\triangleq}

\renewcommand{\Bbb}{\mathbb}

\newcommand{\Fq}{{{\Bbb F}}_{\!q}}

\newcommand{\cod}{{\rm cod}}

\theoremstyle{definition}

\title{Combinatorial alphabet-dependent bounds for insdel codes}
\author{Xiangliang~Kong, Itzhak~Tamo and  Hengjia~Wei%
\thanks{This work was supported in part by the European Research Council (ERC) under Grant 852953 and the National Natural Science Foundation of China under Grant 12371523.}
\thanks{X. Kong (rongxlkong@gmail.com) and I. Tamo (tamo@tauex.tau.ac.il) are with the Department of Electrical Engineering-Systems, Tel Aviv University, Tel Aviv-Yafo 6997801, Israel. H. Wei (hjwei05@gmail.com) is with the Peng Cheng Laboratory, Shenzhen 518000, China, and with the School of Mathematics and Statistics, Xi'an Jiaotong University, Xi'an 710049, China, and also with the Pazhou Laboratory (Huangpu), Guangzhou 510555, China.}
}
\date{today}

\begin{document}

\maketitle

\begin{abstract}
    Error-correcting codes resilient to synchronization errors such as insertions and deletions are known as insdel codes. Due to their important applications in DNA storage and computational biology, insdel codes have recently become a focal point of research in coding theory.
    In this paper, we present several new combinatorial upper and lower bounds on the maximum size of $q$-ary insdel codes. Our main upper bound is a sphere-packing bound obtained by solving a linear programming (LP) problem. It improves upon previous results for cases when the distance $d$ or the alphabet size $q$ is large. Our first lower bound is derived  from a connection between insdel codes and matchings in special hypergraphs. This lower bound, together with our upper bound, shows that for fixed block length $n$ and edit distance $d$, when $q$ is sufficiently large, the maximum size of insdel codes is $ \frac{q^{n-\frac{d}{2}+1}}{{n\choose \frac{d}{2}-1}}(1 \pm o(1))$. The second lower bound refines Alon et al.'s recent logarithmic improvement on Levenshtein's GV-type bound and extends its applicability to large $q$ and $d$. 
\end{abstract}

\section{Introduction}

Insdel codes, designed for synchronization channels to correct errors arising from both insertions and deletions (insdel for short), have attracted significant attentions from coding theorists recently due to their important applications in DNA storage and computational biology \cite{MBT10,BLCCSS16,YGM17,HMG19}. A code $\cC\subseteq [q]^{n}$ is called an \emph{$(n,d)_q$ insdel code with edit distance $d$} if for any two codewords in $\cC$, the smallest number of insertions and deletions needed to transform one codeword into the other is at least $d$. An $(n,d)_q$ insdel code can correct up to $\lfloor\frac{d-1}{2}\rfloor$ insertion/deletion errors. Denote $D_q(n,d)$ as the maximum size of an $(n,d)_q$ insdel code.

The systematic study of bounds on insdel codes was first initiated by Levenshtein \cite{Levenshtein66} for the binary case. Then, building on the construction using VT codes \cite{VT65}, various binary insdel codes have been developed that aim to reduce the gap to Levenshtein's upper bound. For example, see \cite{HF02, APFC11}, and for recent progress, see \cite{BGZ17, SB20, GH21, GHL22}. Non-binary insdel codes were first studied by Calabi and Harnett \cite{CH69} and Tanaka and Kasai \cite{TK76}.  Following the works of \cite{CH69} and \cite{TK76}, Tenengolts \cite{Tenengolts84} proposed a construction of single-insertion/deletion-correcting codes, i.e., $(n,4)_q$ insdel codes, and showed that they are asymptotically optimal. Later, through combinatorial design theory, perfect $(n,d)_q$ insdel codes with large $d$, such as $d=2n-2$, $2n-4$ and $2n-6$, were constructed and studied in \cite{Bours95,Mahmoodi98,Yin01,SWY02,CGL10,WG15}. For general $n$, $q$ and $d$, Levenshtein \cite{Levenshtein02} extended his result in \cite{Levenshtein66} and proved 
a Gilbert-Varshamov-type (GV-type) lower bound and a sphere-packing-type upper bound on $D_q(n,d)$. He also developed asymptotic bounds on $D_q(n,4)$, as $\frac{q}{n}$ approaches infinity \cite{Levenshtein92}. Levenshtein's upper bound in \cite{Levenshtein02} bound was improved by Kulkarni and Kiyavash in \cite{KK13} for $4\le d\le n+1$, and by Fazeli et al. \cite{FVY15} and by Cullina and Kiyavash \cite{CK16} for $d=4$. Following these works, Yasunaga \cite{Yasunaga24} proved an Elias-type upper bound and an MRRW-type upper bound on $D_q(n,d)$, which improves Levenshtein's upper bound in \cite{Levenshtein02} and Kulkarni and Kiyavash's result in \cite{KK13} for the binary and the quaternary cases. Moreover, in a recent work by Liu and Xing \cite{LX23}, the authors proved several bounds on $D_q(n,d)$, which showed that the insdel-metric Singleton bound is not tight for nontrivial codes. Very recently, Alon et al.  \cite{ABGHK24} obtained the first asymptotic improvement to Levenshtein's GV-type bound and showed the existence of $(n,d)_q$ insdel codes logarithmically larger than the GV-type bound in \cite{Levenshtein02} for fixed $q$ and $d$.

Besides  the studies on bounds for insdel codes,  numerous excellent explicit constructions for non-binary insdel codes have been proposed in recent years \cite{LN16,SGB20,HS21,CCGKN21,SC23}. However, the sizes of the insdel codes obtained in most of these constructions are still significantly smaller than the GV-type lower bound given in \cite{Levenshtein02}. Therefore, improvements on current lower and upper bounds on $D_q(n,d)$ are not only theoretically interesting but may also provide insights for constructing better codes.


In the same spirit as Levenshtein's works \cite{Levenshtein02, Levenshtein92}, Kulkarni and Kiyavash's work \cite{KK13}, and Alon et al.'s work \cite{ABGHK24}, this paper investigates the behavior of $D_q(n,d)$ for general values of $n$, $q$, and $d$, and presents several new upper and lower bounds. Our contributions are as follows:


\begin{itemize}
\item {\bf Upper Bound:} Our upper bound is obtained by reducing the sphere-packing problem for  $(n,d)_q$ insdel codes to a linear programming (LP) problem. This bound is valid for $4\leq d \leq 2n-2$. For $d=2n-2\tau$ with $\tau$ fixed, this bound improves upon Levenshtein's sphere-packing-type bound in \cite{Levenshtein02} by a factor of $O(\frac{1}{n^{\tau+1}})$. For fixed $n$ and $d$ and sufficiently large $q$, this bound improves  both Levenshtein's bound in \cite{Levenshtein02} and Kulkarni and Kiyavash's bound in \cite{KK13} by a constant factor; moreover, it can be shown that this bound is  asymptotically optimal in this case.
\item {\bf First Lower Bound:} Our first lower bound is obtained by relating the existence of an $(n,d)_q$ insdel code to the existence of a matching in a special hypergraph. This lower bound asymptotically achieves  our  upper bound for fixed $n$ and $d$ and sufficiently large $q$. Specifically, if $n$ and $d$ are fixed such that $4\leq d\leq 2n-2$ and $q$ is sufficiently large, then 
\[D_q(n,d)=\frac{q^{n-\frac{d}{2}+1}}{{n\choose \frac{d}{2}-1}}(1 \pm o(1)).\]
It is worth noting that this somewhat  generalizes Levenshtein's result in \cite[Corollary~5.1]{Levenshtein92}, which states that $D_q(n,4)\sim q^{n-1}/n$ as $q/n \to \infty$.  
\item {\bf Second Lower Bound:} Our second lower bound refines Alon et al.'s recent logarithmic improvement \cite{ABGHK24} on Levenshtein's GV-type bound by a constant factor and extends its applicability to large $q$ and $d$. Specifically, when $q\geq n$, this bound improves  upon Levenshtein's GV-type bound by a factor of $(\frac{d}{2}-1)\ln{n}$ if  $d$ is fixed, and by a factor of $\parenv{\frac{\delta}{2}(\ln{\frac{1}{\delta}}-8.1)}n$ if $d-2=\delta n$ for some constant $\delta>0$ such that $\ln \frac{1}{\delta} >8.1$. 
\end{itemize}

The rest of the paper is structured as follows. In Section \ref{sec: preliminaries}, we present several preliminary results and our new bounds on 
$D_q(n,d)$. We also compare our bounds with the previous ones.  Section \ref{sec: upper bounds} is dedicated to proving our upper bounds on the size of 
$q$-ary insdel codes. Following that, in Section \ref{sec: lower bounds}, we establish our lower bounds through two existence results of 
$q$-ary insdel codes.

\section{Preliminaries and main results}\label{sec: preliminaries}

\subsection{Notations and preliminary results}

For integers $1\le m\le n$, let $[m,n]\triangleq \{m,m+1,\dots,n\}$ and  $[n]\triangleq [1,n]$. 
Let $q\geq 2$ be an integer and $\mathcal{C}\subseteq [q]^n$ be a code over $[q]$ of length $n$. 
For a codeword $\bbc\in \mathcal{C}$
and a subset $R\subseteq [n]$, let $\bbc|_{R}$ be the vector obtained by projecting the coordinates of $\bbc$ to $R$, and define $\mathcal{C}|_{R}=\{\mathbf{c}|_{R}: \mathbf{c}\in \mathcal{C}\}$. 
When $q$ is a prime power, we use $\mathbb{F}_q$ to denote the finite field of size $q$ and use $\mathbb{F}_q^n$ to denote the $n$-dimensional linear space over $\mathbb{F}_q$. 
Let $\alpha_1,\alpha_2,\ldots,\alpha_n$, where $n\leq q$, be $n$ pairwise distinct elements of $\mathbb{F}_q$. For $k\in [n]$, denote by $\mathbb{F}_q^{<k}[x]$ the set of polynomials of degree less than $k$. The Reed-Solomon code with evaluation vector $\bm{\alpha}\triangleq (\alpha_1,\alpha_2,\ldots,\alpha_n)\in \Fq^n$ and dimension $k$ is defined by
\begin{equation}\label{eq_def_RS-code}
    RS_{\bm{\alpha}}\triangleq \{(f(\alpha_1),f(\alpha_2),\ldots,f(\alpha_n)):f(x)\in \mathbb{F}_q^{<k}[x]\}.
\end{equation}

Throughout the paper, we use standard asymptotic notation, as follows. Let $f(n)$ and $g(n)$ be two non-negative functions defined on the positive integers. We say that $f=O_{n}(g)$ (or $g=\Omega_{n}(f)$) if there is some constant $c\geq 0$ such that $f(n)\le cg(n)$ for all $n\geq 1$ and we say that $f=o_{n}(g)$ if for every $\epsilon>0$ there exists a constant $n_0$ such that $f(n)\le \epsilon g(n)$ for all $n\geq n_0$. Often, the subscript $n$ will be omitted if it's clear from the context. Unless otherwise specified, logarithm functions $\log(\cdot)$ are base-$2$. 

For two words $\bbx=(x_1,x_2,\ldots,x_m)\in [q]^{m}$ and $\bby = (y_1,y_2,\ldots,y_n)\in [q]^{n}$ with $m\leq n$, $\bbx$ is called a \emph{subsequence} of $\bby$ if there are $1\le i_1<i_2<\cdots<i_m\le n$ such that $\bbx=\bby|_{\{i_1,i_2,\ldots,i_m\}}$. The \emph{edit distance} $d_E(\bbx,\bby)$ between $\bbx$ and $\bby$ is the minimum number of insertions and deletions required to transform $\bbx$ into $\bby$. 
A \emph{longest common subsequence} (LCS) of $\bbx$ and $\bby$ is a common subsequence of both $\bbx$ and $\bby$ that achieves the maximum length. We denote by $LCS(\bbx,\bby)$ the length of an LCS of $\bbx$ and $\bby$. Then, 
the edit distance between $\bbx$ and $\bby$ can be calculated via $LCS(\bbx,\bby)$.
\begin{lemma}\label{lem_ed_distance}(\cite{CR02}, Lemma 12.1)
     Let $\bbx\in [q]^{m}$ and $\bby\in [q]^{n}$, then it holds that     \begin{equation}\label{eq_ed_distance}
        d_E(\bbx,\bby)=m+n-2LCS(\bbx,\bby).
    \end{equation}
\end{lemma}

The edit distance of a $q$-ary code $\cC\subseteq [q]^n$ is defined to be $\min\{d_{E}(\bbu,\bbv):\bbu\neq \bbv\in \cC\}$. A code over $[q]$ of length $n$ with edit distance $d$ is called an $(n,d)_q$ \emph{insdel code}. As is shown by Levenshtein in \cite{Levenshtein66}, an insdel code $\cC\subseteq [q]^n$ can correct $t$ deletion errors if and only if it can correct $t$ insertion errors. Due to this equivalence, we use insdel error-correcting capability to denote deletion and insertion error-correcting capability. Clearly, an $(n,d)_q$ insdel code can correct up to $\lfloor\frac{d-1}{2}\rfloor$ insdel errors.

For positive integers $n\geq t\geq 1$ and $\bbu\in [q]^n$, denote 
$$\cS_{D}(\bbu,t)\triangleq\{\bbv\in [q]^{n-t}:~\bbv~\text{is obtained from}~\bbu~\text{by deleting}~t~\text{symbols}\},$$ 
$$\cS_{I}(\bbu,t)\triangleq\{\bbv\in [q]^{n+t}:~\bbv~\text{is obtained from}~\bbu~\text{by inserting}~t~\text{symbols}\}$$ 
as the \emph{deletion sphere} and \emph{insertion sphere} centered at $\bbx$ of radius $t$, respectively. For even number $2\leq d\leq 2n$, denote 
$$\cB_{E}(\bbu,d)\triangleq\{\bbv\in [q]^{n}:~d_{E}(\bbu,\bbv)\leq d\}$$
as the \emph{editing ball} centered at $\bbx$ of radius $t$. Then, we have the following bounds on $|\cS_{D}(\bbu,t)|$, $|\cS_{I}(\bbu,t)|$ and $|\cB_{E}(\bbu,d)|$.
\begin{lemma}\label{bounds_on_deletion/insertion_ball}
Let $n,t$ be positive integers $n\geq t\geq 1$. Then, for any $\bbu\in [q]^n$, it holds that:
    \begin{itemize}
        \item [1)] ${\rho(\bbu)-t+1\choose t}\leq |\cS_{D}(\bbu,t)|\leq {\rho(\bbu)+t-1\choose t}$, where $\rho(\bbu)$ is the number of runs\footnote{A run is a maximum interval of $\bbu$ consisting of the same symbols.} in $\bbu$.
        \item [2)] $|\cS_{I}(\bbu,t)|=\sum_{i=0}^{t}{n+t\choose i}(q-1)^{i}$.
        \item [3)] $|\cB_{E}(\bbu,2t)|\leq \sum_{i=0}^{t}\left(|\cS_{D}(\bbu,i)|\cdot|\cS_{I}(\bbu,i)|\right)\leq n^{2t}q^t(1+o_{n}(1))$ as $n\rightarrow \infty$.
    \end{itemize}
\end{lemma}
\begin{IEEEproof}
For 1) and 2), the reader is referred to \cite{Levenshtein66} and \cite{Levenshtein01}, respectively. For 3), the bound $|\cB_{E}(\bbu,2t)|\leq \sum_{i=0}^{t}|\cS_{D}(\bbu,i)|\cdot|\cS_{I}(\bbu,i)|$ follows directly from the definitions of $\cS_{D}(\bbu,t)$, $\cS_{I}(\bbu,t)$ and $\cB_{E}(\bbu,d)$. By $\rho(\bbu)\leq n$ and 1), we have $$\sum_{i=0}^{t}|\cS_{D}(\bbu,i)|\leq \sum_{i=0}^{t}{n+i-1\choose i}={n+t\choose t}\leq n^t(1+o_n(1)).$$ 
    Then, by $|\cS_{I}(\bbu,i)|\leq |\cS_{I}(\bbu,t)|$, we have $|\cB_{E}(\bbu,2t)|\leq n^t(1+o_n(1))\cdot|\cS_{I}(\bbu,t)|$. Meanwhile, by 2), 
    $$|\cS_{I}(\bbu,t)|\leq q^t\sum_{i=0}^{t}{n+t\choose i}\leq q^tn^t(1+o_n(1)),$$
    where the second inequality follows by $\sum_{i=0}^{t}{n+t\choose i}\leq \frac{n^2}{2}(1+o_n(1))$ when $t\leq 2$ and $\sum_{i=0}^{t}{n+t\choose i}\leq \left(\frac{\mathrm{e}(n+t)}{t}\right)^{t}=n^t\parenv{\frac{e}{t}+\frac{e}{n}}^{t}\leq n^t(1+o_n(1))$ when $t\geq 3$.\footnote{The inequality $\sum_{i=0}^{t}{n+t\choose i}\leq \left(\frac{\mathrm{e}(n+t)}{t}\right)^{t}$ comes from $\sum_{i=0}^{t}{n+t\choose i}\leq \sum_{i=0}^{t}\frac{(n+t)^i}{i!}=\sum_{i=0}^{t}\frac{t^i}{i!}\parenv{\frac{n+t}{t}}^{i}\leq \parenv{\frac{n+t}{t}}^{t}\sum_{i=0}^{t}\frac{t^i}{i!}\leq \left(\frac{\mathrm{e}(n+t)}{t}\right)^{t}$, where the last inequality follows from $e^t=\sum_{i=0}^{\infty}\frac{t^i}{i!}$.} 
    This confirms 3).
\end{IEEEproof}

\subsection{Upper bounds on $D_q(n,d)$}

In \cite{Levenshtein02}, Levenshtein extended his results for binary insdel codes in \cite{Levenshtein66} and obtained the following sphere-packing type upper bound for $q$-ary insdel codes. 

\begin{theorem}\label{thm_upper_bounds_Levenshtein02}(\cite{Levenshtein02}, Theorem 2)
    For any integers $n$, $q\geq 2$, $d$ and $r$ such that $4\leq d\leq 2n-2$ and $\frac{d}{2}-2\leq r\leq n-1$, it holds that \begin{equation}\label{eq0_Levenshtein02_bound}
    D_{q}(n,d)\leq \frac{q^{n-\frac{d}{2}+1}}{\sum_{i=0}^{\frac{d}{2}-1}{r+2-\frac{d}{2}\choose i}}+q\sum_{i=0}^{r-1}{n-1\choose i}(q-1)^{i}.
    \end{equation}
\end{theorem}

During the study of perfect insdel codes through design theory, Bours \cite{Bours95} obtained the following bound on $D_q(n,d)$ that improves Theorem \ref{thm_upper_bounds_Levenshtein02}. Moreover, in \cite{Bours95}, this bound is also shown to be optimal for some infinite families of $q$ when $n=4,5$ (see Theorem 4.8 in \cite{Bours95}).
\begin{theorem}(\cite{Bours95})\label{thm_bours_upper_bound}
For positive integers $n\geq 2$ and $q\geq 2$, it holds that
\begin{equation}\label{eq_Bours'_upper_bound}
    D_{q}(n,2n-2)\leq \lfloor\frac{q}{n}\lfloor\frac{2(q-1)}{n-1}\rfloor\rfloor+q.
\end{equation}
\end{theorem}

In \cite{KK13}, Kulkarni and Kiyavash improved Levenshtein's upper bound in Theorem \ref{thm_upper_bounds_Levenshtein02} when $d\leq n+1$. They modeled the problem of finding the largest $(n,d)_q$ insdel code as a matching problem on a special hypergraph and obtained the following result.
\begin{theorem}\label{KK_upper_bound}(\cite{KK13}, Corollary 4.2)
Let $n$, $d$ and $q$ be positive integers such that $4\leq d\leq n+1$ and $q\geq 2$, it holds that 
    \begin{align}
        D_{q}(n,d)\leq& \sum_{r=3}^{n-\frac{d}{2}+1}\frac{q(q-1)^{r-1}{n-\frac{d}{2}\choose r-1}}{\delta(r,\frac{d}{2}-1)+\sum_{i=d+r-n-3}^{\min\{\frac{d}{2}-3,r-3\}}\delta(r-2,i)} \nonumber\\
        &+q\sum_{r=1}^{2}{n-\frac{d}{2}\choose r-1}(q-1)^{r-1}, \label{eq_KK_upper_bound}
    \end{align}
    where 
    $$\delta(r,\frac{d}{2}-1)=\begin{cases}
    \sum_{i=0}^{\frac{d}{2}-1}{r-\frac{d}{2}+1\choose i}, & r> \frac{d}{2}-1\geq 0,\\
    1, & r=\frac{d}{2}-1\geq 0,\\
    0, & \frac{d}{2}-1<0~\text{or}~\frac{d}{2}-1>r.
    \end{cases}$$
\end{theorem}

In this paper, we reduce the sphere-packing problem of $(n,d)_q$ insdel codes to an LP problem. Then we obtain the following bound on $D_q(n,d)$ by finding feasible solutions to the dual of this LP problem.

\begin{theorem}\label{Thm2_upper_bound}
    Let $n$ and $d$ be positive integers. If  $4\leq d\leq 2n-6$, then we have
    \begin{align}
        D_q(n,d)\leq& \frac{\prod_{i=0}^{n-\frac{d}{2}}(q-i)}{{n\choose \frac{d}{2}-1}}\left(1+\frac{(n-\frac{d}{2}+1)(n-1)}{2(q-n+\frac{d}{2})}\right) \nonumber\\
        &+U_{n,q,d}, \label{eq6_coro2}
    \end{align}
    where $U_{n,q,d}=\sum_{j=\ceilenv{\frac{n-d/2+1}{2}}+1}^{n-d/2-1}\frac{{q\choose j} j^{n-d/2+1}}{{{2j+d-n-2}\choose{2j+d/2-n-1}}}+\sum_{j=1}^{\ceilenv{\frac{n-d/2+1}{2}}}{q\choose j} j^{n-d/2+1}$. For $d=2n-4$ or $2n-2$, we have
    \begin{align}
        D_q(n,2n-2)&\leq q+\frac{2q(q-1)}{n(n-1)},\label{eq3_coro2}\\
        D_q(n,2n-4)&\leq 3q^2-q+\frac{q(q-1)(q-2)}{{n\choose 3}}. \label{eq4_coro2}
    \end{align}
\end{theorem}

Recall that as $r$ increases, the first term $\frac{q^{n-\frac{d}{2}+1}}{\sum_{i=0}^{\frac{d}{2}-1}{r+2-\frac{d}{2}\choose i}}$ of Levenshtein's upper bound in (\ref{eq0_Levenshtein02_bound}) decreases, while the second term $q\sum_{i=0}^{r-1}{n-1\choose i}(q-1)^{i}$ increases. Thus, the RHS of (\ref{eq0_Levenshtein02_bound}) is at least 
\begin{equation}\label{eq4_Levenshtein02_upper_bound}
    \frac{q^{n-\frac{d}{2}+1}}{\sum_{i=0}^{\frac{d}{2}-1}{n+1-\frac{d}{2}\choose i}}+q\sum_{i=0}^{\frac{d}{2}-3}{n-1\choose i}(q-1)^{i}.
\end{equation}
Based on this lower bound of the RHS of (\ref{eq0_Levenshtein02_bound}), in the following, we compare our bounds in Theorem \ref{Thm2_upper_bound} with Levenshtein's Theorem \ref{thm_upper_bounds_Levenshtein02} and Kulkarni and Kiyavash's Theorem \ref{KK_upper_bound}. As a result, one can see that Theorem \ref{Thm2_upper_bound} improves upon results of Theorem \ref{thm_upper_bounds_Levenshtein02} and Theorem \ref{KK_upper_bound} in the parameter regimes where $d$ or $q$ is large. 

\subsubsection*{The case of $d=2n-2\tau$ with $\tau\geq 1$ fixed}
When $\tau=1$, we have $\frac{d}{2}=n-1$ and $n-\frac{d}{2}+1=2$. Then, for $n\geq 4$, (\ref{eq4_Levenshtein02_upper_bound}) is at least 
    $$\frac{q^{2}}{\sum_{i=0}^{n-2}{2\choose i}}+q=\frac{q^2}{4}+q>\frac{2q(q-1)}{n(n-1)}+q=\text{RHS of } (\ref{eq3_coro2}).$$ 
When $\tau=2$, we have $\frac{d}{2}=n-2$ and $n-\frac{d}{2}+1=3$. Then, for $n\geq 6$, (\ref{eq4_Levenshtein02_upper_bound}) is at least
    \begin{align*}
        \frac{q^{3}}{\sum_{i=0}^{n-3}{3\choose i}}+q\parenv{(n-1)(q-1)+1}&=\frac{q^3}{8}+(n-1)q^2-(n-2)q\\
        &> \frac{q(q-1)(q-2)}{{n\choose 3}}+3q^2-q=\text{RHS of } (\ref{eq4_coro2}).
    \end{align*}
When $\tau=3$, we have $\frac{d}{2}=n-3$ and $n-\frac{d}{2}+1=4$. Thus, for $n\geq 8$, (\ref{eq4_Levenshtein02_upper_bound}) is at least 
    \begin{align}
        &\frac{q^{4}}{\sum_{i=0}^{n-4}{4\choose i}}+q\sum_{i=0}^{2}{n-1\choose i}(q-1)^{i}\notag \\
        =&\frac{q^4}{16}+{n-1\choose 2}q(q-1)^2+{n-1\choose 1}q(q-1)+q. \label{eq:tau3levbnd}
    \end{align} 
    While the bound (\ref{eq6_coro2}) by Theorem \ref{Thm2_upper_bound} is
    \begin{align*}
        D_{q}(n,d)\leq &\frac{\prod_{i=0}^{3}(q-i)}{{n\choose 4}}\left(1+\frac{2(n-1)}{q-3}\right)+\sum_{j=1}^{2}{q\choose j}j^{4}\\
        \leq &\frac{q^4}{{n\choose 4}}+\frac{2(n-1)q(q-1)^2}{{n\choose 4}}+8q(q-1)+q,
    \end{align*} 
which is strictly smaller than~\eqref{eq:tau3levbnd}. In general, for $d=2n-2\tau$, the upper bounds in Theorem \ref{Thm2_upper_bound} improves Levenshtein's upper bound (\ref{eq0_Levenshtein02_bound}) by a factor of $O(\frac{1}{n^{\tau+1}})$.

\subsubsection*{The case where $n$ and $d$ are fixed and $q\rightarrow\infty$}
{Since $n$ and $d$ are fixed, we have that $j^{n-d/2+1}=\Theta(1)$ and ${{2j+d-n-2}\choose{2j+d/2-n-1}}=\Theta(1)$ for $1\leq j\leq n-\frac{d}{2}-1$. This leads to $U_{n,q,d}=O(q^{n-\frac{d}{2}-1})$. Then, by (\ref{eq6_coro2}), we have 
$$D_q(n,d) \leq \frac{q^{n-\frac{d}{2}+1}}{{n\choose \frac{d}{2}-1}}\parenv{1+\Theta\parenv{\frac{1}{q}}}+O(q^{n-\frac{d}{2}-1})$$
for $4\leq d\leq 2n-6$. Combining with (\ref{eq3_coro2}) and (\ref{eq4_coro2}), it follows that
\begin{equation}\label{eq_upper_bound_q_infty}
    D_q(n,d)\leq\frac{q^{n-\frac{d}{2}+1}}{{n\choose \frac{d}{2}-1}}(1+o(1)),
\end{equation}
when $n$ and $d$ are fixed and $q\rightarrow\infty$.
}

Note that when $r=n-\frac{d}{2}+1$, we have $d+r-n-3=\frac{d}{2}-2>\min\{\frac{d}{2}-3,r-3\}$. Thus, $\sum_{i=d+r-n-3}^{\min\{\frac{d}{2}-3,r-3\}}\delta(r-2,i)=0$ when $r=n-\frac{d}{2}+1$. Then, for fixed $n$ and $d$ such that $n\geq4$, as $q\rightarrow\infty$, the upper bound (\ref{eq_KK_upper_bound}) in Theorem \ref{KK_upper_bound} becomes
\begin{equation}\label{eq_KK's_SP_upper_bound}
D_{q}(n,d)\leq \frac{q^{n-\frac{d}{2}+1}}{\delta(n-\frac{d}{2}+1,\frac{d}{2}-1)}(1+o_q(1)).
\end{equation} 

Consider the lower bound (\ref{eq4_Levenshtein02_upper_bound}) on the RHS of Levenshtein's upper bound in (\ref{eq0_Levenshtein02_bound}). When $4\leq d\leq n+2$, we have $\frac{d}{2}-2<n-\frac{d}{2}+1$. Then, the first term in (\ref{eq4_Levenshtein02_upper_bound}) is the dominating term. Thus, in this case, (\ref{eq4_Levenshtein02_upper_bound}) becomes
\begin{equation}\label{eq1_Levenshtein02_upper_bound}
    \frac{q^{n-\frac{d}{2}+1}}{\sum_{i=0}^{\frac{d}{2}-1}{n-d+2\choose i}}(1+o_q(1)).
\end{equation}
Note that $\delta(n-\frac{d}{2}+1,\frac{d}{2}-1)=\sum_{i=0}^{\frac{d}{2}-1}{n-d+2\choose i}$, thus (\ref{eq_KK's_SP_upper_bound}) reduces to Levenshtein's bound (\ref{eq1_Levenshtein02_upper_bound}). 
When $n+3\leq d\leq 2n-2$, we have $\frac{d}{2}-2\geq n-\frac{d}{2}+1$. Then, the second term in (\ref{eq4_Levenshtein02_upper_bound}) is the dominating term. Thus, in this case, (\ref{eq4_Levenshtein02_upper_bound}) is lower bounded by
\begin{equation}\label{eq2_Levenshtein02_upper_bound}
    {n-1\choose {\frac{d}{2}-3}}q^{\frac{d}{2}-2}(1-o_q(1)).
\end{equation}

Therefore, by 
$$\sum_{i=0}^{\frac{d}{2}-1}{n-d+2\choose i}<\sum_{i=0}^{\frac{d}{2}-1} {n-d+2 \choose i} {d-2 \choose \frac{d}{2}-1-i} = {n \choose \frac{d}{2}-1},$$ 
the bound (\ref{eq_upper_bound_q_infty}) of Theorem \ref{Thm2_upper_bound} improves upon both Levenshtein's \cite{Levenshtein02} and Kulkarni and Kiyavash's results \cite{KK13} when $n,d$ are fixed and $q\rightarrow\infty$. This improvement is not marginal. For example, when $d=n+1$, we have $\sum_{i=0}^{\frac{d}{2}-1}{n-d+2\choose i}=\sum_{i=0}^{\frac{d}{2}-1}{1\choose i}=2$ and (\ref{eq1_Levenshtein02_upper_bound}) becomes $q^{n-\frac{d}{2}+1}(\frac{1}{2}+o_q(1))$. By (\ref{eq_upper_bound_q_infty}), Theorem \ref{Thm2_upper_bound}  
 provides a constant factor improvement. When $d\geq n+3$, Kulkarni and Kiyavash's result Theorem \ref{KK_upper_bound} no longer holds and Levenshtein's upper bound in Theorem \ref{thm_upper_bounds_Levenshtein02} is at least (\ref{eq2_Levenshtein02_upper_bound}), which is even larger than (\ref{eq1_Levenshtein02_upper_bound}).

\vspace{10pt}

It is also worth noting that since $D_{q}(n,2n-2)$ is an integer, (\ref{eq3_coro2}) in Theorem \ref{Thm2_upper_bound} implies that $D_{q}(n,2n-2)\leq \lfloor\frac{2q(q-1)}{n(n-1)}\rfloor+q$. This coincides with the Bours' upper bound in Theorem \ref{thm_bours_upper_bound} asymptotically. Moreover, as an application of the upper bounds  on $D_q(n,2n-2)$ and $D_q(n,2n-4)$, in Section \ref{sec: upper bound 2}, we also obtain a $d$-dependent factor improvement to the second bound in the following result by Liu and Xing in \cite{LX23} for $4\leq d\leq \min\{2q-2,2n-2\}$. 
\begin{theorem}\label{thm_Liu_Xing}(\cite{LX23}, Theorem 3.6)
    For integers $d$ and $n$ such that $2\leq d\leq 2n$, the following bounds hold for $D_q(n,d)$. 
    \begin{itemize}
        \item [1.] $D_q(n,2)=q^n$ and $D_q(n,2n)=q$.
        \item [2.] $D_q(n,d)\leq \frac{1}{2}(q^{n-\frac{d}{2}+1}+q^{n-\frac{d}{2}})$ for $4\leq d\leq 2n-2$.
        \item [3.] $D_q(n,d)\leq q^{n-\frac{d}{2}}$ for $2q\leq d\leq 2n-2$.
    \end{itemize}
\end{theorem}

\subsection{Lower bounds on $D_q(n,d)$}

Levenshtein \cite{Levenshtein02} proved the following GV-type lower bound for $q$-ary insdel codes. 

\begin{theorem}\label{thm_lower_bounds_Levenshtein02}(\cite{Levenshtein02}, Theorem 1)
    For any integers $n$, $q\geq 2$ and $d$ such that $4\leq d\leq 2n-2$, it holds that \begin{equation}\label{eq1_Levenshtein02_bound}
    D_{q}(n,d)\geq \frac{q^{n+\frac{d}{2}-1}}{(\sum_{i=0}^{\frac{d}{2}-1}{n\choose i}(q-1)^{i})^{2}}.
    \end{equation}
\end{theorem}

This lower bound is known as the best lower bound on $D_q(n,d)$ for general $n$, $q$ and $d$.
Recently,  Alon et al. \cite{ABGHK24} proved the following improved lower bound on $D_q(n,d)$ for fixed alphabet size $q$ and fixed distance $d$. 

\begin{theorem}\label{log_improvement_Alon}(\cite{ABGHK24}, Theorem 1)
 For fixed $d$ and $q$ such that $6\leq d\leq 2n+2$ and $q\geq 2$, it holds that
    $$D_q(n,d)\geq \Omega_{q,d}(q^n\log{n}/n^{d-2}).$$
\end{theorem} 

Note that when both $q$ and $d$ are fixed, Levenshtein's lower bound in Theorem \ref{thm_lower_bounds_Levenshtein02} reduces to
\begin{equation}\label{eq1_Levenshtein02_lower_bound}
D_{q}(n,d)\geq \frac{q^{n+\frac{d}{2}-1}}{n^{d-2}(q-1)^{d-2}}(1-o_n(1)).
\end{equation} 
Thus, in this case, Alon et al.'s result in \cite{ABGHK24} is logarithmically larger than the bound in Theorem \ref{thm_lower_bounds_Levenshtein02}.

In this paper, we present two new lower bounds on $D_q(n,d)$. The first lower bound is obtained by reducing the existence of an $(n,d)_q$ insdel code to the existence of a matching in a special hypergraph.
\begin{theorem}\label{Thm_optimla_codes_existence}
    Let $n$ and $d$ be fixed positive integers such that $4\leq d\leq 2n-2$. Then, we have 
    \begin{equation}\label{eq_lower_bound_q_infty}
        D_{q}(n,d)\geq \frac{q^{n-\frac{d}{2}+1}}{{n\choose \frac{d}{2}-1}}(1-o(1)),
    \end{equation}
    as $q\rightarrow\infty$.
\end{theorem}
When $q$ is sufficiently large, Levenshtein's lower bound in Theorem \ref{thm_lower_bounds_Levenshtein02} is at most
\begin{equation}\label{eq_Levenshtein02_lower_bound}
    \frac{q^{n-\frac{d}{2}+1}}{{n\choose \frac{d}{2}-1}^2}(1-o_q(1)).
\end{equation}
Compared to Levenshtein's lower bound in (\ref{eq_Levenshtein02_lower_bound}), when both $n$ and $d$ are fixed, the lower bound (\ref{eq_lower_bound_q_infty}) in Theorem \ref{Thm_optimla_codes_existence} provides a ${n\choose \frac{d}{2}-1}$ factor improvement. Moreover, by \eqref{eq_upper_bound_q_infty}, this bound is asymptotically  optimal. 
\begin{corollary}\label{cor:fix n and d}
Let $n$ and $d$ be fixed positive integers such that $4\leq d\leq 2n-2$. Then, we have 
\[D_q(n,d)=  \frac{q^{n-\frac{d}{2}+1}}{{n\choose \frac{d}{2}-1}}(1 \pm o_q(1)).\]
as $q\rightarrow\infty$.
\end{corollary}

The second lower bound is a refinement of Alon et al.'s Theorem \ref{log_improvement_Alon} in \cite{ABGHK24}.

\begin{theorem}\label{Thm_GV_type_lower_bound}
    Let $n$, $q\geq 2$ and $d$ be positive integers.
    \begin{itemize}
        \item If $q\geq n$ is prime power and $6\leq d< \frac{n}{\mathrm{e}^{8.1}}+2$, then it holds that
        \begin{equation}\label{eq_GV_type_lower_bound}
        D_q(n,d)\geq \parenv{\frac{d}{2}-1}(\ln{n}-\ln{(d-2)}-8.1)(1-o(1))\frac{q^{n-\frac{d}{2}+1}}{{n\choose \frac{d}{2}-1}^2},
        \end{equation}
         as $n\rightarrow \infty$.
        \item If $q\geq 2$ and $d\geq 6$ are fixed, then it holds that
        \begin{equation}\label{eq_GV_type_lower_bound2}
        D_q(n,d)\geq \parenv{\frac{d}{2}-1-o_n(1)}\frac{q^{n-\frac{d}{2}+1}\log{n}}{n^{d-2}}, 
        \end{equation}
        as $n\rightarrow \infty$.
    \end{itemize}
\end{theorem}
When  $d$ is fixed and $q$ increases with $n$, the lower bound (\ref{eq_GV_type_lower_bound}) in Theorem \ref{Thm_GV_type_lower_bound} offers a $(\frac{d}{2}-1)\ln{n}$ factor improvement over (\ref{eq_Levenshtein02_lower_bound}). When $d-2=\delta n$ for some constant $\delta>0$ such that $\ln \frac{1}{\delta} >8.1$, (\ref{eq_GV_type_lower_bound}) provides a $\parenv{\frac{\delta}{2}(\ln{\frac{1}{\delta}}-8.1)}n$ factor improvement over (\ref{eq_Levenshtein02_lower_bound}). 

When both $q$ and $d$ are fixed, the lower bound (\ref{eq_GV_type_lower_bound2}) in Theorem \ref{Thm_GV_type_lower_bound} refines the bound in Theorem \ref{log_improvement_Alon} by specifying the constant factor. 

\section{Upper bounds on the size of insdel codes}\label{sec: upper bounds}

In this section, we prove our upper bounds on the size of the maximum size $D_q(n,d)$ of an $(n,d)_q$ insdel code. As an application, we obtain an result that improves Liu and Xing's Singleton-type bound in \cite{LX23}. 

\subsection{Proof of Theorem \ref{Thm2_upper_bound}}\label{sec: upper bound 1}

Let $4\leq d \leq 2n-2$ be an even number and $\cC\subseteq [q]^{n}$ be an $(n,d)_q$ insdel code. By the definition of insdel code, for any word $\bbu\in [q]^{n-\frac{d}{2}+1}$, there is at most one codeword in $\cC$ containing $\bbu$ as a subsequence. That is, for any two distinct words $\bbc,\bbc'\in \cC$, it holds that
\begin{equation*}\label{eq00_upper_bound2}
    \cS_{D}(\bbc,\frac{d}{2}-1)\cap \cS_{D}(\bbc',\frac{d}{2}-1)=\emptyset.
\end{equation*}
Denote $s\triangleq n-\frac{d}{2}$, we have $1\leq s \leq n-2$ as $4\leq d \leq 2n-2$. This implies that each $(n,d)_q$ insdel code induces a sphere-packing of $[q]^{s+1}$. Based on this observation, we prove the first upper bound by reducing this packing problem to a linear programming problem.

\begin{theorem}\label{Thm_upper_bound}
    Let $q$, $n$ and $s$ be positive integers such that $1\leq s\leq n-2$. Let $M$ be the maximum value of the following LP problem:
    \begin{align}
        \textrm{\rm maximize}\quad & \sum_{i=1}^{n}x_i \nonumber\\
        \textrm{\rm subject to} \quad \sum_{i=1}^{\lceil\frac{s+1}{2}\rceil-1} & x_i\leq \sum_{j=1}^{\lceil\frac{s+1}{2}\rceil-1} b_j; \label{eq2_Thm_upper_bound2}\\
        \sum_{i=j}^{n+j-(s+1)}a_{i,j} & x_i\leq b_{j}, ~\forall~ \lceil\frac{s+1}{2}\rceil\leq j\leq s+1;\label{eq1_Thm_upper_bound2}\\
        & x_i\geq 0,~\forall~ 1\leq i\leq n,\nonumber
    \end{align}
    where\footnote{For positive integer $n$ and $k\in \mathbb{Z}$, we set ${n\choose k}=0$ if $k<0$ or $k>n$. For nonegative integer $a_1,a_2,\ldots,a_i$, ${n\choose a_1,a_2,\ldots,a_i}$ is defined as the number $\frac{n!}{a_1!\cdots a_i!(n-a_1-\cdots-a_i)!}$.} 
    $a_{i,j}\triangleq {{i+j-(s+1)}\choose {2j-(s+1)}}$ and $b_j\triangleq \sum_{c_1,\ldots,c_j\in \mathbb{Z}_{+}\atop c_1+\cdots+c_{j}=s+1}{q\choose j}{s+1\choose c_1,c_2,\ldots,c_j}.$
    Then, we have $D_q(n,2n-2s)\leq M$. 
\end{theorem}

Before delving into the detailed proof of Theorem \ref{Thm_upper_bound}, we first prove Theorem \ref{Thm2_upper_bound} by finding feasible solutions to the dual LP problem corresponding to the LP problem stated in Theorem \ref{Thm_upper_bound}.

\begin{IEEEproof}[Proof of Theorem \ref{Thm2_upper_bound}]
    The dual LP problem of the LP problem in Theorem \ref{Thm_upper_bound} is:
    \begin{align}
        \text{minimize}\quad & \left(\sum_{j=1}^{\lceil\frac{s+1}{2}\rceil-1} b_j\right)y_0+\sum_{j=\ceilenv{\frac{s+1}{2}}}^{s+1}b_{j} y_j\nonumber\\
        \text{subject to} \quad &y_0\geq 1~\text{when}~s\geq 2,\nonumber\\
        \sum_{j=\ceilenv{\frac{s+1}{2}}}^{\min\{i,s+1\}}a_{i,j}&y_j\geq 1,~\forall~\ceilenv{\frac{s+1}{2}}\leq i\leq n-\floorenv{\frac{s+1}{2}},\label{eq1_coro2}\\
        \sum_{j=i+s+1-n}^{\min\{i,s+1\}}a_{i,j}&y_j\geq 1,~\forall~n-\floorenv{\frac{s+1}{2}}+1\leq i\leq n,\label{eq2_coro2}\\
        &y_i\geq 0,~\forall~i=0,\ceilenv{\frac{s+1}{2}}\leq i\leq s+1.\nonumber
    \end{align}

    When $s=1$, we have $\lfloor\frac{s+1}{2}\rfloor=\lceil\frac{s+1}{2}\rceil=1$. Then, one can easily verify that $y_0=0$, $y_1=a_{1,1}^{-1}=1$ and $y_2=a_{n,2}^{-1}={n\choose 2}^{-1}$ is a feasible solution of the dual LP problem. Thus, by the weak duality theorem (see Proposition 6.1.1 in \cite{Matouvsek2007}), $\sum_{i=1}^{n}x_i$ is upper bounded by 
    $$b_{1}y_1+b_2y_2=q+\frac{2q(q-1)}{n(n-1)}.$$
    This confirms (\ref{eq3_coro2}).

    When $s=2$, we have $\lfloor\frac{s+1}{2}\rfloor=1$ and $\lceil\frac{s+1}{2}\rceil=2$. Then, one can easily verify that $y_0=1$, $y_2=a_{2,2}^{-1}=1$ and $y_3=a_{n,3}^{-1}={n\choose 3}^{-1}$ is a feasible solution of the dual LP problem. Thus, 
    $$\sum_{i=1}^{n}x_i\leq b_{1}y_0+b_2y_2+b_3y_3=q+{q\choose 2} 3!+\frac{{q\choose 3}}{{n\choose 3}} 3!.$$
    This confirms (\ref{eq4_coro2}).
    
    For general $s\geq 3$, $y_0=1$, $y_{\ceilenv{\frac{s+1}{2}}}=a_{\ceilenv{\frac{s+1}{2}},\ceilenv{\frac{s+1}{2}}}^{-1}=1$ and
    \begin{equation}\label{eq5_coro2}
        \begin{cases}
            y_{s+1}=a_{n,s+1}^{-1},\\
            y_{s}=a_{n-1,s}^{-1}(1-\frac{a_{n-1,s+1}}{a_{n,s+1}}),\\
            y_{j}=a_{n+j-(s+1),j}^{-1},~\ceilenv{\frac{s+1}{2}}+1\leq j\leq s-1
        \end{cases}
    \end{equation} 
    is a feasible solution of the dual LP problem: $y_0=1$ satisfies the first condition, $y_{\ceilenv{\frac{s+1}{2}}}=a_{\ceilenv{\frac{s+1}{2}},\ceilenv{\frac{s+1}{2}}}^{-1}=1$ implies that (\ref{eq1_coro2}) holds naturally, and one can then easily verify that (\ref{eq2_coro2}) holds by (\ref{eq5_coro2}).
    Noting that $b_{s+1}={q\choose s+1} (s+1)!$, $b_s={q\choose s} \frac{s(s+1)!}{2}$ and by $\sum_{c_1,\ldots,c_j\in \mathbb{Z}_{\geq 0}\atop c_1+\cdots+c_{j}=s+1}{s+1\choose c_1,c_2,\ldots,c_j}=j^{s+1}$, 
    \begin{equation}\label{eq_bound_b_j}
        b_j=\sum_{c_1,\ldots,c_j\in \mathbb{Z}_{+}\atop c_1+\cdots+c_{j}=s+1}{q\choose j}{s+1\choose c_1,c_2,\ldots,c_j}\leq {q\choose j} j^{s+1}.
    \end{equation}
    Then, (\ref{eq5_coro2}) implies that $\sum_{i=1}^{n}x_i$ is at most
    \begin{align*} 
    &\left(\sum_{j=1}^{\lceil\frac{s+1}{2}\rceil-1} b_j\right)y_0+\sum_{j=\ceilenv{\frac{s+1}{2}}}^{s+1}b_{j} y_j \nonumber\\ 
    \overset{(I)}{=}& \frac{b_{s+1}}{a_{n,s+1}}+\frac{b_{s}}{a_{n-1,s}}(1-\frac{a_{n-1,s+1}}{a_{n,s+1}})+\sum_{j=\ceilenv{\frac{s+1}{2}}+1}^{s-1}\frac{b_{j}}{a_{n+j-(s+1),j}}+\sum_{j=1}^{\lceil\frac{s+1}{2}\rceil} b_j \nonumber\\
    \overset{(II)}{\leq}& \frac{  {q\choose s+1} (s+1)! }{{n\choose s+1}}\left(1+\frac{(s+1)(n-1)}{2(q-s)}\right)+\sum_{j=\ceilenv{\frac{s+1}{2}}+1}^{s-1}\frac{{q\choose j} j^{s+1}}{{{n+2j-2(s+1)}\choose{2j-(s+1)}}}+\sum_{j=1}^{\ceilenv{\frac{s+1}{2}}}{q\choose j} j^{s+1},
    \end{align*}
    where $(I)$ follows by $y_0=y_{\ceilenv{\frac{s+1}{2}}}=1$ and (\ref{eq5_coro2}), and (II) follows by $a_{i,j}= {{i+j-(s+1)}\choose {2j-(s+1)}}$ and (\ref{eq_bound_b_j}). This confirms (\ref{eq6_coro2}).
\end{IEEEproof}

    

Next, for the proof of Theorem \ref{Thm_upper_bound}, we introduce some extra notations and preliminaries results.

For positive integer $n$, let $\bbu=(u_1,u_2,\ldots,u_n)$ be a word in $[q]^n$. We denote $N(\bbu)\triangleq \{\alpha\in [q]:\exists~i\in [n]~\text{such that}~u_i=\alpha\}$ as the set of symbols appear in $\bbu$. For each $\alpha\in N(\bbu)$, we denote $F_{\bbu}(\alpha)\triangleq|\{i\in [n]: u_i=\alpha\}|$ as the frequency that $\alpha$ appears in $\bbu$.
Assume that $N(\bbu)=\{\alpha_1,\alpha_2,\ldots,\alpha_{|N(\bbu)|}\}$ and $F_{\bbu}(\alpha_1)\geq \cdots\geq F_{\bbu}(\alpha_{|N(\bbu)|})$. Clearly, we have $F_{\bbu}(\alpha_{|N(\bbu)|})\geq 1$ and $\sum_{i=1}^{|N(\bbu)|}F_{\bbu}(\alpha_i)=n$. In the following, for $1\leq i\leq n-1$, we abbreviate $\cS_{D}(\bbu,n-i)$ as $\cS_{i}(\bbu)$. That is, we use $\cS_{i}(\bbu)$ to denote the set of all length $i$ subsequences of $\bbu$ in $[q]^{i}$. 


\begin{lemma}\label{lem_upper_bound2}
    Let $a$, $b$ and $c$ be positive integers such that $b\leq a\leq 2b$, $b\leq c$ and $a-b+c\leq n$. Let $\bbu=(u_1,u_2,\ldots,u_n)\in[q]^n$ with $N(\bbu)=\{\alpha_1,\alpha_2,\ldots,\alpha_{c}\}$. Then, for any $(2b-a)$-subset $\{\beta_{1},\ldots,\beta_{2b-a}\}\subseteq\{\alpha_{a-b+1},\ldots,\alpha_{c}\}$, there is a word $\bby\in \cS_a(\bbu)$ with $N(\bby)=\{\beta_{1},\ldots,\beta_{2b-a}\}\cup \{\alpha_{1},\ldots,\alpha_{a-b}\}$ and $F_{\bby}(\beta_1)=\cdots=F_{\bby}(\beta_{2b-a})=1$.
\end{lemma}

\begin{IEEEproof}
    Recall that $F_{\bbu}(\alpha_1)\geq \cdots\geq F_{\bbu}(\alpha_{c})$. Assume that $\sum_{i=1}^{a-b}F_{\bbu}(\alpha_i)=n_1$, then, we claim that $n_1\geq 2(a-b)$. Otherwise, by $n_1\leq 2(a-b)-1$ and $F_{\bbu}(\alpha_i)\geq 1$, we must have $1=F_{\bbu}(\alpha_{a-b})=F_{\bbu}(\alpha_{a-b+1})=\cdots=F_{\bbu}(\alpha_c)$. This implies that $\sum_{i=a-b+1}^{c}F_{\bbu}(\alpha_i)=c-a+b=n-n_1$. Hence, we have 
    \begin{align*}
        n&\leq c-a+b+2(a-b)-1\\
        &= a-b+c-1,
    \end{align*}
    which contradicts the last condition on $a$, $b$ and $c$.

    Let $I\subseteq [n]$ be the set of indices such that $u_{i}\in \{\alpha_{a-b+1},\ldots,\alpha_{c}\}$ for every $i\in I$. Then, we have $|I|=\sum_{i=a-b+1}^{c}F_{\bbu}(\alpha_i)$, and by $F_{\bbu}(\alpha_i)\geq 1$, we know that $|I|\geq c-a+b \geq 2b-a$. Pick $2b-a$ indices $i_1,i_2,\ldots,i_{2b-a}$ from $I$ such that $u_{i_j}=\beta_j$ for every $1\leq j\leq 2b-a$. Note that $|[n]\backslash I|=n_1\geq 2(a-b)$. Thus, we can pick a subset $I'$ of $2(a-b)$ indices from $[n]\backslash I$ such that $N(\bbu|_{I'})=\{\alpha_{1},\ldots,\alpha_{a-b}\}$. Then, let $I'' \triangleq \{i_1,i_2,\ldots,i_{2b-a}\} \cup I'$. Clearly, $|I''|=a$ and it is easy to check that $\bby \triangleq \bbu|_{I''}$ is the desired word.
\end{IEEEproof}

With the help of Lemma \ref{lem_upper_bound2}, we can proceed the proof of Theorem \ref{Thm_upper_bound}.

\begin{IEEEproof}[Proof of Theorem \ref{Thm_upper_bound}]
    Let $\cC\subseteq [q]^{n}$ be an $(n,2n-2s)_q$ insdel code. Then, for any word $\bbu\in [q]^{s+1}$, there is at most one codeword in $\cC$ containing $\bbu$ as a subsequence. That is, for distinct $\bbc,\bbc'\in \cC$, it holds that
    \begin{equation}\label{eq0_upper_bound2}
        \cS_{s+1}(\bbc)\cap \cS_{s+1}(\bbc')=\emptyset.
    \end{equation}

    According to the size of $N(\bbu)$ for each $\bbu\in [q]^{s+1}$, we can divide $[q]^{s+1}$ into $s+1$ parts, $[q]^{s+1}=\bigcup_{j=1}^{s+1}\cA_{j}$, where $\cA_{j}=\{\bbu\in [q]^{s+1}: |N(\bbu)|=j\}$. Then, we have $|\cA_{j}|=b_j.$
    Similarly, we can divide $\cC$ into $n$ parts, $\cC=\bigcup_{i=1}^{n}\cC_i$, where $\cC_i=\{\bbc\in\cC: |N(\bbc)|=i\}$. For each $i\in [n]$, we denote $\cS_{j}(\cC_i)=\bigcup_{\bbc\in\cC_i}\cS_{j}(\bbc)$ as the set of all length $j$ subseqences of codewords in $\cC_{i}$. Thus, by (\ref{eq0_upper_bound2}), we have $|\cC_{i}|\leq |\cS_{s+1}(\cC_i)|$. Next, we proceed the proof by showing that $(|\cC_1|,\ldots,|\cC_n|)$ lies in the feasible region of the stated LP problem.
    
    First, we show that 
    \begin{equation}\label{eq1_upper_bound2}
        \cS_{s+1}(\cC_i)\subseteq 
        \begin{cases}
        \bigcup_{j=s+1-(n-i)}^{\min\{s+1,i\}}\cA_{j},~\text{if}~s+1> n-i,\\
        \bigcup_{j=1}^{\min\{s+1,i\}}\cA_{j},~\text{otherwise}.
        \end{cases}
    \end{equation}
    Since $\cS_{s+1}(\cC_i)\subseteq [q]^{s+1}$ and $N(\bbc)=i$ for each $\bbc\in \cC_i$, it holds naturally that $\cS_{s+1}(\cC_i)\subseteq \bigcup_{j=1}^{\min\{s+1,i\}}\cA_{j}$ for every $1\leq i\leq n$. When $s+1> n-i$, we claim that any length $s+1$ subsequence $\bbu$ of some codeword $\bbc\in \cC_i$ satisfies $N(\bbu)\geq s+1-(n-i)$. Otherwise, we would have $N(\bbc)\leq s-(n-i)+n-(s+1)=i-1$, which contradicts the definition of $\cC_i$. Thus, in this case, we have $\cS_{s+1}(\cC_i)\subseteq\bigcup_{j=s+1-(n-i)}^{\min\{s+1,i\}}\cA_{j}$.
    
    For $1\leq i\leq \lceil\frac{s+1}{2}\rceil-1$, by (\ref{eq1_upper_bound2}) and $\lceil\frac{s+1}{2}\rceil-1\leq s-1$, we have
    $$\cS_{s+1}(\cC_i)\subseteq \bigcup_{j=1}^{i}\cA_{j}\subseteq \bigcup_{j=1}^{\lceil\frac{s+1}{2}\rceil-1}\cA_{j}.$$
    Recall that $|\cC_{i}|\leq|\cS_{s+1}(\cC_i)|$. This leads to $\sum_{i=1}^{\lceil\frac{s+1}{2}\rceil-1}|\cC_{i}|\leq \sum_{j=1}^{\lceil\frac{s+1}{2}\rceil-1}b_j$.
    Therefore, $|\cC_1|,\ldots,|\cC_n|$ satisfy condition (\ref{eq2_Thm_upper_bound2}).
    
    Noticed that if $\cS_{s+1}(\cC_i)\cap \cA_{j}\neq \emptyset$, then by (\ref{eq1_upper_bound2}), we have
    \begin{equation*}
        \begin{cases}
        s+1-(n-i)\leq j\leq \min\{s+1,i\},~\text{when}~i>n-(s+1),\\
        1\leq j\leq \min\{s+1,i\},~\text{when}~1\leq i\leq n-(s+1).
        \end{cases}
    \end{equation*}
    That is, $\cS_{s+1}(\cC_i)\cap \cA_{j}\neq \emptyset$ only if $i\leq n+j-(s+1)$. Next, we show that for each $i\in [n]$, $j\in [s+1]$ satisfying $2j\geq s+1$ and $j\leq i\leq n+j-(s+1)$, it holds that
    \begin{equation}\label{eq3_upper_bound2}
        |\cS_{s+1}(\cC_i)\cap \cA_{j}|\geq |\cC_{i}|{{i+j-(s+1)}\choose {2j-(s+1)}}.
    \end{equation}
    Let $\bbc\in\cC_i$ and assume that $N(\bbc)=\{\alpha_1,\alpha_2,\ldots,\alpha_{i}\}$. Applying Lemma \ref{lem_upper_bound2} with $\bbu=\bbc$, $a=s+1$, $b=j$ and $c=i$. Then, for each $2j-(s+1)$-subset $\{\beta_{1},\ldots,\beta_{2j-(s+1)}\}\subseteq\{\alpha_{s-j+2},\ldots,\alpha_{i}\}$, there is a word $\bby\in \cS_{s+1}(\bbc)$ with $N(\bby)=\{\beta_{1},\ldots,\beta_{2j-(s+1)}\}\cup \{\alpha_{1},\ldots,\alpha_{s+1-j}\}$ and $F_{\bby}(\beta_1)=\cdots=F_{\bby}(\beta_{2j-(s+1)})=1$. Clearly, $\bby\in \cS_{s+1}(\bbc)\cap \cA_{j}$. Note that for different choices of $\{\beta_{1},\ldots,\beta_{2j-(s+1)}\}\subseteq\{\alpha_{s-j+2},\ldots,\alpha_{i}\}$, the corresponding $\bby\in \cS_{s+1}(\bbc)\cap \cA_{j}$ are also different. Thus, we have
    \begin{equation}\label{eq4_upper_bound2}
        |\cS_{s+1}(\bbc)\cap \cA_{j}|\geq {{i+j-(s+1)}\choose {2j-(s+1)}}.
    \end{equation}
    Meanwhile, by (\ref{eq0_upper_bound2}), we have $|\cS_{s+1}(\cC_i)\cap \cA_{j}|=\sum_{\bbc\in\cC_i}|\cS_{s+1}(\bbc)\cap \cA_{j}|$. Therefore, (\ref{eq3_upper_bound2}) follows directly from (\ref{eq4_upper_bound2}).

    For $\lceil\frac{s+1}{2}\rceil \leq i\leq n$ and any fixed $j\in [s+1]$, (\ref{eq3_upper_bound2}) implies that
    $$\sum_{i=j}^{n+j-(s+1)}|\cC_{i}| {{i+j-(s+1)}\choose {2j-(s+1)}}\leq \sum_{i=j}^{n+j-(s+1)}|\cS_{s+1}(\cC_i)\cap \cA_{j}|.$$
    Since $\cC_1,\ldots,\cC_n$ is a partition of $\cC$, by (\ref{eq0_upper_bound2}), we know that $\cS_{s+1}(\cC_1),\ldots,\cS_{s+1}(\cC_n)$ are also pairwise disjoint. Thus, for every $\lceil\frac{s+1}{2}\rceil\leq j\leq s+1$, it holds that
    \begin{equation}\label{eq5_upper_bound2}
    \sum_{i=j}^{n+j-(s+1)}|\cC_{i}| {{i+j-(s+1)}\choose {2j-(s+1)}}\leq |b_{j}|.
    \end{equation}
    This shows that $|\cC_1|,\ldots,|\cC_n|$ satisfy condition (\ref{eq1_Thm_upper_bound2}).
    
    Finally, the conclusion follows directly by $|\cC_i|\geq 0$ for every $1\leq i\leq n$. 
\end{IEEEproof}



\subsection{An application of Theorem \ref{Thm2_upper_bound}}\label{sec: upper bound 2}

In this section, we first prove a recursive bound for $D_q(n,d)$. Then, as an application of this recursive bound and the upper bounds (\ref{eq_Bours'_upper_bound}) on $D_q(n,2n-2)$ and (\ref{eq4_coro2}) on $D_q(n,2n-4)$, we obtain another upper bound on $D_q(n,d)$ that improves the result of Theorem \ref{thm_Liu_Xing} for $4\leq d\leq \min\{2q-2,2n-4\}$.

\begin{proposition}\label{prop_recursive_bound}
    Let $n$, $q\geq 2$ and $d$ be positive integers such that $4\leq d\leq 2n-2$ and $d$ is even. For integer $1\leq i\leq n-\frac{d}{2}$, suppose that $D_q(\frac{d}{2}+i,d)\leq f(d,q,i)$ for some function of $d$, $q$ and $i$ satisfying $f(d,q,i)\geq q^i$ for any $d$. Then, we have
    $$D_q(n,d)\leq q^{n-\frac{d}{2}-i}f(n,q,i).$$
\end{proposition}

\begin{IEEEproof}
Let $\cC\subseteq [q]^{n}$ be an insdel code with edit distance $d$ and Hamming distance $d'$. Suppose that 
\begin{equation}\label{eq1_prop_recursive_bound}
   |\cC|\geq q^{n-\frac{d}{2}-i}f(n,q,i)+1. 
\end{equation}
Next, we proceed with the proof by deriving a contradiction.

First, we show that $d=2d'$. Otherwise, since $d\leq 2d'$ (see Lemma 3.1 in \cite{LX23}), this implies that $d\leq 2d'-2$. Thus, we have $q^{n-\frac{d}{2}-i}\geq q^{n-d'+1-i}$. By (\ref{eq1_prop_recursive_bound}) and $f(n,q,i)\geq q^{i}$, this further leads to $|\cC|\geq q^{n-d'+1}+1$, which contradicts to the Singleton bound for codes under the Hamming metric. Thus, we have $d=2d'$.

For each $i\in [q]$, define set
$$\cB_{i}\triangleq\{(c_1,\ldots,c_n)\in \cC: c_1=i\}.$$
Clearly, $\cB_1,\cB_2,\ldots,\cB_q$ are pairwise disjoint and form a partition of $\cC$. Thus, we have $|\cC|=\sum_{i=1}^{q}|\cB_i|\leq q\cdot \max_{1\leq i\leq q}|\cB_i|$. Let $l\in [q]$ such that  $|\cB_l|= \max_{1\leq i\leq q}|\cB_i|$. Then we have $|\cB_l|\geq \frac{|\cC|}{q}$. 

Let $\cC_1$ be the code obtained from $\cB_l$ by deleting the first coordinate, i.e., $\cC_1=\cB_l|_{[2,n]}$. Then, $\cC_1$ is a code of length $n-1$ with Hamming distance $d_{H}(\cC_1)\geq d'$ and size at least $\frac{|\cC|}{q}$. Next, we claim that $\cC_1$ has edit distance $d_E(\cC_1)\geq d$. Let $\bbu,\bbv\in [q]^{n-1}$ be two distinct codewords of $\cC_1$ satisfying 
$$d_E(\cC_1)=2(n-1)-2LCS(\bbu,\bbv).$$
By the definition of $\cB_l$, $(l,\bbu)$ and $(l,\bbv)$ are two distinct codewords of $\cC$. Thus, we have 
\begin{align*}
    d&\leq 2n-2LCS((l,\bbu),(l,\bbv))\\
     &= 2n-2(1+LCS(\bbu,\bbv))=d_E(\cC_1).
\end{align*}

Continuing this fashion, we can obtain a sequence of codes $\{\cC_{i}\}_{i=0}^{n-(d'+i)}$ with $\cC_0=\cC$, and for $i=0,1,\ldots,n-(d'+i)$, it holds that: 1) $\cC_i\subseteq [q]^{n-i}$; 2) $|\cC_i|\geq \frac{|\cC|}{q^{i}}$; 3) $d_{H}(\cC_i)\geq d'$ and $d_E(\cC_i)\geq d$.
Specially, the code $\cC_{n-(d'+i)}$ is a code of length $d'+i$ with Hamming distance $d_{H}(\cC_{n-(d'+i)})\geq d'$ and the edit distance $d_{E}(\cC_{n-(d'+i)})\geq d$. Moreover, the size of $\cC_{n-(d'+i)}$ is at least 
\begin{align}\label{eq2_prop_recursive_bound}
    |\cC_{n-(d'+i)}|&\geq \frac{|\cC|}{q^{n-(d'+i)}} \geq f(d,q,i)+\frac{1}{q^{n-(d'+i)}},
\end{align}
where the last inequality follows from (\ref{eq1_prop_recursive_bound}) and $d'=\frac{d}{2}$.

Next, we show that $d_{H}(\cC_{n-(d'+i)})= d'$ and $d_{E}(\cC_{n-(d'+i)})= d$. Suppose $d_{H}(\cC_{n-(d'+i)})\geq d'+1$, then, by the Singleton bound for codes under Hamming metric, $|\cC_{n-(d'+i)}|\leq q^{d'+i-d_{H}(\cC_{n-(d'+i)})+1}\leq q^i$. This contradicts to (\ref{eq2_prop_recursive_bound}). Thus, $d_{H}(\cC_{n-(d'+i)})= d'$. Then, this leads to $d_{E}(\cC_{n-(d'+i)})\leq 2d_{H}(\cC_{n-(d'+i)})= 2d'=d$. On the other hand, we know that $d_{E}(\cC_{n-(d'+i)})\geq d$. Thus, $d_{E}(\cC_{n-(d'+i)})= d$. 

Finally, the contradiction follows directly by (\ref{eq2_prop_recursive_bound}) and $|\cC_{n-(d'+i)}|\leq D_q(\frac{d}{2}+i,d)\leq f(d,q,i)$.  
\end{IEEEproof}

Applying Proposition \ref{prop_recursive_bound} with $i=1$ and the bound on $D_q(n,2n-2)$ in (\ref{eq_Bours'_upper_bound}), and with $i=2$ and the bound on $D_q(n,2n-4)$ in (\ref{eq4_coro2}), we can obtain the following upper bound on $D_q(n,d)$.

\begin{theorem}\label{thm_improved_LX_upper_bound}
    Let $n$, $q\geq 2$ and $d$ be positive integers such that $4\leq d\leq 2n-4$ and $d$ is even. Then, it holds that
    $$D_{q}(n,d)\leq q^{n-\frac{d}{2}-1}\cdot\min\left\{\lfloor\frac{2q}{d+2}\lfloor\frac{4(q-1)}{d}\rfloor\rfloor+q,\frac{48(q-1)(q-2)}{(d+4)(d+2)d}+3q-1\right\}.$$
\end{theorem}



\begin{remark}\label{rmk3}

Note that for $d\geq 4$,
\begin{align*}
    q^{n-\frac{d}{2}-1}\parenv{\lfloor\frac{2q}{d+2}\lfloor\frac{4(q-1)}{d}\rfloor\rfloor+q}&\leq 
     q^{n-\frac{d}{2}-1}\parenv{\frac{8q(q-1)}{(d+2)d}+q }\\
    &< q^{n-\frac{d}{2}-1}{\frac{q^2+q}{2}}
\end{align*}
Thus, Theorem~\ref{thm_improved_LX_upper_bound} improves Liu and Xing's bound $D_q(n,d)\leq \frac{1}{2}(q^{n-\frac{d}{2}+1}+q^{n-\frac{d}{2}})$. 

We also note that when $d\geq 2q$, $\lfloor\frac{2q}{d+2}\lfloor\frac{4(q-1)}{d}\rfloor\rfloor+q=q$ and $\frac{48(q-1)(q-2)}{(d+4)(d+2)d}+3q-1\geq 3q-1$, thus the bound in Theorem~\ref{thm_improved_LX_upper_bound} reduce to
$D_{q}(n,d)\leq q^{n-\frac{d}{2}-1}\cdot q=q^{n-\frac{d}{2}}$,
which coincides with Liu and Xing's bound $D_q(n,d)\leq q^{n-\frac{d}{2}}$ in Theorem \ref{thm_Liu_Xing}.

Moreover, in a recent work by Liu et al. \cite{LWGZ24}, the authors also studied the behavior of $D_q(n,d)$ and obtained several new upper bounds on $D_q(n,d)$ for $n \geq q$. Specifically, they showed that
$$D_q(n,d)\leq q^{n-\frac{d}{2}-\frac{d}{2q-2}}(q+2)$$
when $q \mid n$ and $2q-2 \leq d \leq 2n-\frac{2n}{q}$. This improves Theorem \ref{thm_Liu_Xing} by Liu and Xing and Theorem \ref{thm_improved_LX_upper_bound} for the case when $q\mid n$, $2q \leq d \leq 2n-2$ and $q^{\frac{d}{2q-2}} \geq q+2$.
\end{remark}

\section{Lower bounds on the size of Insdel codes}\label{sec: lower bounds}

In this section, we present the proofs of our lower bounds on the size of the maximum size $D_q(n,d)$ of an $(n,d)_q$ insdel code.

\subsection{Proof of Theorem \ref{Thm_optimla_codes_existence}}\label{sec: lower bound1}

Very recently, using Kahn's theorem \cite{Kahn96} on the existence of near-optimal matchings, Liu and Shangguan \cite{LSG24} proved the existence of near-optimal constant weight codes and constant composition codes for all fixed odd distances. In the following,  building on a hypergraph characterization of insdel codes given by Kulkarni and Kiyavash \cite{KK13}, we use  a similar idea to show the existence of optimal insdel codes (w.r.t. the bound in Theorem \ref{Thm2_upper_bound}) for all fixed $n$ and $4\leq d\leq 2n-2$ when $q$ is sufficient large.

Recall that a hypergraph $\cH$ is a tuple $(V,\cE)$, where the vertex set $V$ is a finite set and the hyperedge set $\cE$ is a collection of nonempty subsets of $V$. We denote $\Delta(\cH)$ as the \emph{maximum degree} of $\cH$ and $\text{cod}(\cH)$  as the \emph{maximum codegree} of $\cH$, i.e., $\Delta(\cH)\eqdef\max_{v\in V}|\{E\in \cE: v\in E\}|$ and $\cod(\cH)\eqdef\max_{v_1\neq v_2\in V}|\{E\in \cE: v_1,v_2\in E\}|$. We call $\cH$ an $l$-bounded hypergraph if $|E|\leq l$ for every $E\in \cE$. A matching of a hypergraph $\cH=(V,\cE)$ is a collection of pairwise disjoint hyperedges $E_1,\ldots,E_j\in \cE$. The matching number of $\cH$, denoted as $\nu(\cH)$, is the largest $j$ for which such a matching exists.

Consider the following hypergraph:
\begin{equation}\label{eq1_insdel_hypergraph}
    \cH_{q,t,n}^{D}\triangleq([q]^{n-t},\{\cS_D(\bbu,t):~\bbu\in [q]^{n}\}).
\end{equation}
In $\cH_{q,t,n}^{D}$, vertices are all words in $[q]^{n-t}$ while hyperedges correspond to length-$(n-t)$ subsequences of a word $\bbu$ in $[q]^{n}$. As pointed out in \cite{KK13}, an insdel code that can correct $t$ deletions  corresponds to a matching in $\cH_{q,t,n}^{D}$. Recall that an $(n,d)_q$ insdel code can correct at most $\frac{d}{2}-1$ deletions/insertions. Thus, we have $D_q(n,d)=\nu(\cH_{q,\frac{d}{2}-1,n}^{D})$.

To obtain a lower bound on $\nu(\cH_{q,\frac{d}{2}-1,n}^{D})$, we need some extra results. The first is the following variation of Kahn's theorem on the existence of near-optimal matchings in \cite{LSG24}.

\begin{lemma}\label{lem_nearly-optimal_matchings}(\cite{LSG24},~Lemma 2.4)
    Let $\cH$ be an $l$-bounded hypergraph with $\frac{\cod(\cH)}{\Delta(\cH)}=o(1)$ as $|V(\cH)|\rightarrow\infty$. Then, we have
    $$\nu(\cH)>(1-o(1))\frac{|\cE(\cH)|}{\Delta(\cH)}.$$
\end{lemma}

In order to employ Lemma \ref{lem_nearly-optimal_matchings}, we also need the following bounds on the maximum degree and maximum codegree of $\cH_{q,t,n}^{D}$.

\begin{lemma}\label{lem_degree_codegree_insdel_hypergraph}
For positive integer $n$, $t$ and $q$ such that $t\leq n-1$, the following holds for $\cH_{q,t,n}^{D}$:
\begin{itemize}
    \item [1)] $\cH_{q,t,n}^{D}$ is ${n\choose t}$-bounded. 
    \item [2)] $\Delta(\cH_{q,t,n}^{D})=\sum_{i=0}^{t}{n\choose i}(q-1)^{i}$.
    \item [3)]
    $\cod(\cH_{q,t,n}^{D})=\sum_{i=0}^{t-1}{n\choose i}(q-1)^{i}(1-(-1)^{t-i}).$
\end{itemize}
\end{lemma}

\begin{IEEEproof}
Since the size of an hyperedge in $\cH_{q,t,n}^{D}$ is the number of the length-$(n-t)$ subseqences of the corresponding word in $[q]^{n}$, 1) holds naturally. Moreover, by the definition of $\cH_{q,t,n}^{D}$, we know that
\begin{align*}
    \Delta(\cH_{q,t,n}^{D})&=\max_{\bbv\in [q]^{n-t}}|\cS_{I}(\bbv,t)|,\\
    \cod(\cH)&=\max_{\bbv\neq \bbv'\in [q]^{n-t}}|\cS_{I}(\bbv,t)\cap \cS_{I}(\bbv',t)|,
\end{align*}
where $\cS_{I}(\bbv,t)$ is the insertion sphere centered at $\bbv$ of radius $t$. Therefore, 2) and 3) follows directly from $|\cS_{I}(\bbv,t)|=\sum_{i=0}^{t}{n\choose i}(q-1)^{i}$ by Lemma \ref{bounds_on_deletion/insertion_ball} (see also equation (24) in \cite{Levenshtein01}), and 
$$\max_{\bbv\neq \bbv'\in [q]^{n-t}}|\cS_{I}(\bbv,t)\cap \cS_{I}(\bbv',t)|=\sum_{i=0}^{t-1}{n\choose i}(q-1)^{i}(1-(-1)^{t-i})$$ 
by Theorem 3 and equation (51) in \cite{Levenshtein01}.
\end{IEEEproof}

The proof of Theorem \ref{Thm_optimla_codes_existence} is now immediate.

\begin{IEEEproof}[Proof of Theorem \ref{Thm_optimla_codes_existence}]
By Lemma \ref{lem_degree_codegree_insdel_hypergraph},
$$\frac{\cod(\cH_{q,t,n}^{D})}{\Delta(\cH_{q,t,n}^{D})}=\frac{\sum_{i=0}^{t-1}{n\choose i}(q-1)^{i}(1-(-1)^{t-i})}{\sum_{i=0}^{t}{n\choose i}(q-1)^{i}}=\Theta\parenv{\frac{1}{q}}\rightarrow 0,$$
as $q\rightarrow\infty$. Applying Lemma \ref{lem_nearly-optimal_matchings} on $\cH_{q,\frac{d}{2}-1,n}^{D}$, we have
\begin{align*}
    \nu(\cH)>(1-o(1))\frac{|\cE(\cH)|}{\Delta(\cH)}\geq \frac{q^{n-\frac{d}{2}+1}}{{n\choose \frac{d}{2}-1}}(1-o(1)),
\end{align*}
as $q\rightarrow\infty$.
\end{IEEEproof}

\subsection{Proof of Theorem \ref{Thm_GV_type_lower_bound}}\label{sec: lower bound2}

Alon et al. \cite{ABGHK24} reduced the construction of an $(n,d)_q$-insdel code to identifying an independent set in the deletion graph $\Gamma$. Here, $\Gamma$ comprises vertices from $[q]^n$, with two words connected if their edit distance is at most $d-2$. The key tool employed in \cite{ABGHK24} is the following result by Ajtai et al. \cite{AKZ80} concerning independence numbers of graphs with few triangles.

\begin{lemma}\label{lem_less_triangle_large_independent_set}(\cite{AKZ80}, Lemma 7)
For any $\varepsilon>0$ and any graph $G$ on $N\geq 1$ vertices with average degree $d$ containing $T<Nd^{2-\varepsilon}$ triangles, the independence number 
$$\alpha(G)\geq c_{\varepsilon}\cdot \frac{N\log(d)}{d},$$
where $c_{\varepsilon}= 0.01\cdot \frac{\varepsilon}{48}.$
\end{lemma}

Note that the degree of a vertex $\bbv\in[q]^n$ in $\Gamma$ equals $|\cB_{E}(\bbv,d-2)|$. As pointed in the proof of Theorem \ref{thm_lower_bounds_Levenshtein02} in \cite{Levenshtein02}, the average degree of $\Gamma$, i.e., the average size of $|\cB_{E}(\bbv,d-2)|$, is bounded from above by 
\begin{equation}\label{eq_average_degree_Gamma}
    \frac{(\sum_{i=1}^{\frac{d}{2}-1}{n\choose i}(q-1)^{i})^{2}}{q^{\frac{d}{2}-1}}\leq {n\choose \frac{d}{2}-1}^2q^{\frac{d}{2}-1}(1+o_{n}(1))
\end{equation}
for fixed $d$. Thus, to apply Lemma \ref{lem_less_triangle_large_independent_set}, one needs to estimate the number of triangles in $\Gamma$. In \cite{ABGHK24}, the authors obtained the following bound on the number of triples $(\bbu,\bbv,\bbw)\in ([q]^n)^{3}$ with prescribed values of $d_{E}(\bbu,\bbv)$, $d_{E}(\bbu,\bbw)$ and $d_{E}(\bbv,\bbw)$.

\begin{lemma}\label{triangles_in_Gamma_Alon_orin}(\cite{ABGHK24}, Lemma 7)
    Let $n\geq a\geq b\geq c\geq 1$. The number of triples $(\bbu,\bbv,\bbw)\in ([q]^n)^{3}$ with $d_{E}(\bbu,\bbv)\leq a$, $d_{E}(\bbu,\bbw)\leq b$ and $d_{E}(\bbv,\bbw)\leq c$ is $C_{a,q}\cdot (q^n n^{a+b+c}(\log n)^{b+c-a})$, where $C_{a,q}$ is a constant dependent on $a$ and $q$.
\end{lemma}

When both the alphabet size $q$ and the distance $d$ are fixed, setting $a=b=c=\frac{d}{2}-1$ in Lemma \ref{triangles_in_Gamma_Alon_orin} yields an upper bound of $O_{d,q}(q^nn^{\frac{3d-6}{2}}(\log{n})^{d-2})$ on the number of triangles in $\Gamma$. Then, considering $\abs{V(\Gamma)}=q^n$ and (\ref{eq_average_degree_Gamma}), one can verify that the condition in Lemma \ref{lem_less_triangle_large_independent_set} holds for $\Gamma$ with $\varepsilon=1/3$ as $n\rightarrow \infty$. This results in the logarithmic improvement stated in Theorem \ref{log_improvement_Alon}.

However, when either $q$ or $d$ scales with $n$, the above approach proposed by Alon et al. \cite{ABGHK24} to derive Theorem~\ref{log_improvement_Alon} is inadequate. Specifically, the upper bound on the number of triangles in $\Gamma$ obtained by Lemma \ref{triangles_in_Gamma_Alon_orin} is too large in such cases, which fails to meet the condition in Lemma \ref{lem_less_triangle_large_independent_set}. In the following, we first address this flaw in Lemma \ref{triangles_in_Gamma_Alon_orin}. Then, we propose an approach to circumvent this problem and refine the bound in Theorem \ref{log_improvement_Alon}. 

Applying the same method as in the proof of Lemma~\ref{triangles_in_Gamma_Alon_orin}, we can determine the hidden coefficient in $O_{a,q}(\cdot)$. For detailed calculations, please refer to the appendix.  
\begin{lemma}\label{triangles_in_Gamma_Alon}
    Let $n\geq a\geq b\geq c\geq 1$. The number of triples $(\bbu,\bbv,\bbw)\in ([q]^n)^{3}$ with $d_{E}(\bbu,\bbv)\leq a$, $d_{E}(\bbu,\bbw)\leq b$ and $d_{E}(\bbv,\bbw)\leq c$ is at most
    \begin{align*}
        &n^{a+b+c}q^{n+2b+2c}(66a)^{2b+2c+1}(\log{n})^{2b+2c}(1+o_{n}(1))\\
        &+q^{n+2a}n^{-4a}(3+o_{n}(1)).
    \end{align*}
\end{lemma}

When either $q$ or $d$ is large, the upper bound in Lemma \ref{triangles_in_Gamma_Alon} becomes weak, failing to satisfy the condition outlined in Lemma \ref{lem_less_triangle_large_independent_set}. Specifically, set $a=b=c=\frac{d}{2}-1$, when $q\geq \sqrt{n}$, the first term of the bound in Lemma \ref{triangles_in_Gamma_Alon} is at least
\begin{align*}
    n^{\frac{3d-6}{2}}q^{n+2d-4}\geq q^n\left(n^{2d-4}q^{d-2}\right);
\end{align*}
when $d\geq n^{1/4}+2$, it is at least
\begin{align*}
    n^{\frac{3d-6}{2}}q^{n+2d-4}(33(d-2))^{2d-3}>& \ n^{\frac{3d-6}{2}}q^{n+2d-4}(d-2)^{2d-3}\\
    \geq & \ q^n\left(n^{\frac{8d-15}{4}}q^{2d-4}\right)\\
    > & \ q^n\left(n^{2d-4}q^{d-2}\right).
\end{align*}
Meanwhile, by $\abs{V(\Gamma)}=q^n$ and (\ref{eq_average_degree_Gamma}), the condition in Lemma \ref{lem_less_triangle_large_independent_set} on the number of triangles in $\Gamma$ is $$T\leq q^n\parenv{{n\choose \frac{d}{2}-1}^{2}q^{d-2}}^{2-\varepsilon}<q^n\left(n^{2d-4}q^{d-2}\right).$$ Thus, for both cases, the condition in Lemma \ref{lem_less_triangle_large_independent_set} doesn't hold. 

Similar to Alon et al.'s approach in \cite{ABGHK24}, we also reduce the problem to finding a large independent set of a specific graph. The difference is that the graph we consider is the subgraph of $\Gamma$ induced by some RS-code over $\Fq^{n}$. This graph possesses additional internal algebraic structures that allow us to achieve a better lower bound on the independence number. Moreover, instead of using Lemma \ref{lem_less_triangle_large_independent_set}, we employee the following result by Hurley and Pirot in \cite{PH21}. 

\begin{lemma}[\cite{PH21}, Corollary 1]\label{lem_chi_bounded_neighborhood}
    Let $G$ be a graph of maximum degree $\Delta$ such that for every $x\in V(G)$, the subgraph of $G$ induced by $N_{G}(x)$\footnote{For a simple graph, the notation $N_G(x)=\{y\in V(G): x\sim y~\text{in}~G\}$ denotes the (open) neighborhood of vertex $x$ in $G$, which doesn't include $x$ itself.} has at most $\Delta^2/f$ edges for some $1\leq f\leq \Delta^2+1$. Then, the chromatic number of $G$ satisfies
    $$\chi(G)\leq \left(1+o(1)\right)\frac{\Delta}{\ln(\min\{f,\Delta\})},$$
    as $f\rightarrow \infty$ (and therefore also $\Delta\rightarrow \infty$).
\end{lemma}

As an immediate application of Lemma \ref{lem_chi_bounded_neighborhood}, one can easily obtain a lower bound on the independence number of the graph $G$ defined in Lemma \ref{lem_chi_bounded_neighborhood}:
\begin{equation}\label{eq_lb_independent_bounded_neighborhood}
    \alpha(G)\geq \frac{|V(G)|}{\Delta}\cdot \ln(\min\{f,\Delta\})\cdot \left(1-o(1)\right).
\end{equation}

To introduce the subgraph that will be studied, we require some notions from \cite{ABGHK24}. For $1\leq \lambda\leq n$, we say that $\bbu\in [q]^{n}$ is \emph{$\lambda$-nonrepeating} if $\bbu|_{I}\neq \bbu|_{J}$ for all pairs of distinct intervals $I,J\subseteq [n]$ of length $\lambda$; $\bbu$ is \emph{$\lambda$-repeating} otherwise. Let $n\geq k\geq 1$ be positive integers and $\bm{\alpha}\in\Fq^n$ be a fixed evaluation vector, denote $\Gamma_{k}$ as the subgraph of $\Gamma$ induced by codewords of the RS-code $RS_{\bm{\alpha}}(n,k)$. For any $1\leq \lambda\leq \frac{k-1}{2}$, denote $\Gamma_{k,\lambda}$ as the subgraph of $\Gamma_{k}$ induced by all the $\lambda$-nonrepeating words of $RS_{\bm{\alpha}}(n,k)$. 

Let $k=n-\frac{d}{2}+1$. In the following, we are going to  estimate:
\begin{itemize}
    \item [1.] the number of vertices in $\Gamma_{n-\frac{d}{2}+1,\lambda}$, i.e., the number of $\lambda$-nonrepeating words in $RS_{\bm{\alpha}}(n,n-\frac{d}{2}+1)$;
    \item [2.] the maximum degree of $\Gamma_{n-\frac{d}{2}+1,\lambda}$;
    \item [3.] the number of edges within the neighborhood of any vertex $\bbu\in V(\Gamma_{n-\frac{d}{2}+1,\lambda})$, i.e., the number of pairs $(\bbv,\bbw)$ such that both $\bbv$ and $\bbw$ are $\lambda$-nonrepeating words in $RS_{\bm{\alpha}}(n,n-\frac{d}{2}+1)$ and all $d_E(\bbu,\bbv)$, $d_E(\bbu,\bbw)$ and $d_E(\bbv,\bbw)$ are at most $d-2$.
\end{itemize}
Then we can show that the condition in Lemma \ref{lem_chi_bounded_neighborhood} is fulfilled \footnote{We note that for the subgraph $\Gamma_{n-\frac{d}{2}+1,\lambda}$ with large $q$ and fixed $d$, the condition in Lemma~\ref{lem_less_triangle_large_independent_set} is also satisfied, whereas  Lemma~\ref{lem_chi_bounded_neighborhood} can yield a better bound by a constant factor. } and obtain a  lower bound on the independence number of  $\Gamma_{n-\frac{d}{2}+1,\lambda}$.


\begin{lemma}\label{lem_bound_on_lambda-repaeating_words}
    For $1\leq \lambda\leq \frac{k-1}{2}$, the number of $\lambda$-repeating words in $RS_{\bm{\alpha}}(n,k)$ is at most ${n\choose 2} q^{k-\lambda}$. Hence, the number of vertices of  $\Gamma_{k,\lambda}$ is at least $q^k -{n\choose 2} q^{k-\lambda}$.
\end{lemma}

\begin{IEEEproof}
    Let $f=\sum_{i=0}^{k-1}f_ix^i\in \Fq^{<k}[x]$ such that $
    \bbf=(f(\alpha_1),f(\alpha_2),\ldots,f(\alpha_n))$ is $\lambda$-repeating. Then, there are two distinct intervals $I,J\subseteq [n]$ of length $\lambda$ such that
    $\bbf|_{I}=\bbf|_{J}$.
    W.l.o.g., assume that $I=[i,i+\lambda-1]$ and $J=[j,j+\lambda-1]$ for some $1\leq i< j\leq n+1-\lambda$. Then, we have
    \begin{equation}\label{eq2_lambda-repaeating_words}
        (f_0,f_1,\ldots,f_{k-1})\left(\begin{array}{cccc}
            1 & 1 & \cdots & 1 \\
            \alpha_{i} & \alpha_{i+1} & \cdots & \alpha_{i+\lambda-1} \\
            \vdots & \vdots &  & \vdots \\
            \alpha_{i}^{k-1} & \alpha_{i+1}^{k-1} & \cdots & \alpha_{i+\lambda-1}^{k-1}
        \end{array}\right)=
        (f_0,f_1,\ldots,f_{k-1})\left(\begin{array}{cccc}
             1 & 1 & \cdots & 1 \\
            \alpha_{j} & \alpha_{j+1} & \cdots & \alpha_{j+\lambda-1} \\
            \vdots & \vdots &  & \vdots \\
            \alpha_{j}^{k-1} & \alpha_{j+1}^{k-1} & \cdots & \alpha_{j+\lambda-1}^{k-1}
        \end{array}\right).
    \end{equation}
     Denote  
    $$\bbA\eqdef \left(\begin{array}{cccc}
            \alpha_{i}-\alpha_{j} & \alpha_{i+1}-\alpha_{j+1} & \cdots & \alpha_{i+\lambda-1}-\alpha_{j+\lambda-1} \\
            \alpha^2_{i}-\alpha^2_{j} & \alpha^2_{i+1}-\alpha^2_{j+1} & \cdots & \alpha^2_{i+\lambda-1}-\alpha^2_{j+\lambda-1} \\
            \vdots & \vdots &  & \vdots \\
            \alpha_{i}^{k-1}-\alpha_{j}^{k-1} & \alpha_{i+1}^{k-1}-\alpha_{j+1}^{k-1} & \cdots & \alpha_{i+\lambda-1}^{k-1}-\alpha_{j+\lambda-1}^{k-1}
        \end{array}\right).$$
    Thus, we have $(f_1,\ldots,f_{k-1})\in Ker(\bbA)$, where 
    $Ker(\bbA)\eqdef \{ \bbx \in \Fq^{k-1}: \bbx A =\mathbf{0}\}.$ 
    Note that $Ker(\bbA)$ has size $q^{k-1-rank(\bbA)}$. Thus, given two distinct intervals $I,J\subseteq [n]$ of length $\lambda$, the number of $f\in \Fq^{<k}[x]$ satisfying (\ref{eq2_lambda-repaeating_words}) is at most $q^{k-rank(\bbA)}$. Since the number of distinct intervals $I,J\subseteq [n]$ of length $\lambda$ is at most ${n\choose 2}$, we only need to show that $rank(\bbA)\geq \lambda$.

    Assume that $j=i+l$ for some $l\geq 1$. For $1\leq s\leq n$, denote $\bm{\alpha}_s=(\alpha_s,\alpha_s^2,\ldots,\alpha_s^{k-1})^{T}$. Let $\bbB$ be the $\min\{2\lambda, l+\lambda\}\times \lambda$ matrix of the following form
    $$\bbB=\left(\begin{array}{cccc}
            1 &  &  &  \\
            \vdots & 1 &  &  \\
            -1 & \vdots & \ddots &  \\
             & -1 &  &  \\
             &  & \ddots & 1\\
             &  &  & \vdots \\
             &  &  & -1 
        \end{array}\right),$$
    where the empty parts of $\bbB$ are all zeros and there are $\min \{l-1,\lambda-1\}$ zeros between $1$ and $-1$ in each column. Then, we have $\bbA=\bbA' \bbB$, where 
    $$\bbA'=(\bm{\alpha}_{i},\ldots,\bm{\alpha}_{i-1+\min\{l,\lambda\}},\bm{\alpha}_{j},\ldots,\bm{\alpha}_{j+\lambda-1}).$$ 
    Clearly, $\bbB$ has full column rank. Since $\bbA'$ is the $(k-1)\times \min\{2\lambda, l+\lambda\}$ submatrix obtained by deleting the first row of the $k\times \min\{2\lambda, l+\lambda\}$ Vandermonde matrix generated by $\{{\alpha}_{i},\ldots,{\alpha}_{i-1+\min\{l,\lambda\}},{\alpha}_{j},\ldots,{\alpha}_{j+\lambda-1}\}$, by $\lambda\leq \frac{k-1}{2}$, we have $$rank(\bbA')=\min\{2\lambda,l+\lambda\}.$$
    Thus, by Sylvester's rank inequality (see Theorem 6.5.5 in \cite{Hohn13}), we have 
    \begin{align*}
    rank(\bbA)&\geq rank(\bbA')+rank(B)-\min\{2\lambda,l+\lambda\}=\lambda.
    \end{align*}
    This completes the proof.
\end{IEEEproof}

\begin{lemma}
    The maximum degree of  $\Gamma_{n-\frac{d}{2}+1,\lambda}$ satisfies 
    \begin{equation}\label{eq_maximum_degree_Gamma_{k,lambda}}
    \Delta\parenv{\Gamma_{n-\frac{d}{2}+1,\lambda}} < {n \choose n-\frac{d}{2}+1}^2.
\end{equation}
\end{lemma}

\begin{IEEEproof}
    Since $\Gamma_{n-\frac{d}{2}+1,\lambda}$ is a subgraph of $\Gamma_{n-\frac{d}{2}+1}$, we have $\Delta\parenv{\Gamma_{n-\frac{d}{2}+1,\lambda}}\leq \Delta\parenv{\Gamma_{n-\frac{d}{2}+1}}$. By the definition of $\Gamma_{n-\frac{d}{2}+1}$, for a vertex $\bbu\in RS_{\bm{\alpha}}(n,n-\frac{d}{2}+1)$ of $\Gamma_{n-\frac{d}{2}+1}$, it's neighbor in $\Gamma_{n-\frac{d}{2}+1}$ are the words in $RS_{\bm{\alpha}}(n,n-\frac{d}{2}+1)$ that have edit distance at most $d-2$ with $\bbu$. That is,
    $$\abs{N_{\Gamma_{n-\frac{d}{2}+1}}(\bbu)}=\abs{\left\{\bbv\in RS_{\bm{\alpha}}(n,n-\frac{d}{2}+1):~LCS(\bbu,\bbv)\geq n-\frac{d}{2}+1\right\}}.$$
    On the other hand, $\bbu$ contains at most $n \choose n-\frac{d}{2}+1$ subsequences of length $n-\frac{d}{2}+1$. For each such subsequence $\bbs\in [q]^{n-\frac{d}{2}+1}$ and each subset $I \subset [n]$ of indices with $|I|=n-\frac{d}{2}+1$, there is exactly one codeword $\bbv\in RS_{\bm{\alpha}}(n,n-\frac{d}{2}+1)$ such that $\bbv|_I=\bbs$. Thus, for any $\lambda$, the maximum degree of $\Gamma_{n-\frac{d}{2}+1,\lambda}$ satisfies \begin{equation*}\Delta\parenv{\Gamma_{n-\frac{d}{2}+1,\lambda}}\leq \Delta\parenv{\Gamma_{n-\frac{d}{2}+1}} < {n \choose n-\frac{d}{2}+1}^2.
\end{equation*}
\end{IEEEproof}
We note that the  bound~\eqref{eq_maximum_degree_Gamma_{k,lambda}} for the subgraph $\Gamma_{n-\frac{d}{2}+1,\lambda}$ is independent of the alphabet size $q$, whereas the bound~(\ref{eq_average_degree_Gamma}) for the graph $\Gamma$ depends on $q$.

To establish an upper bound, independent of $q$, on the number of edges in the neighborhood of any vertex $\bbu$ in $\Gamma_{n-\frac{d}{2}+1,\lambda}$,  further concepts from \cite{ABGHK24} are needed. 
We say an element $i$ of a set $I\subseteq [0,n]$ is \emph{$\lambda$-isolated} if $\lambda<i<n-\lambda$ and no other element $j\in I$ satisfies $|j-i|\leq 2\lambda$. Given $\bbu=(u_1,u_2,\ldots,u_n)\in [q]^{n}$, $t\in \{del,ins_1,ins_2,\ldots,ins_q\}$ and $i\in [0,n]$\footnote{Here, $i$ is not allowed to be $0$ when $t=del$.}, we define $\phi_{i,t}(\bbu)$ as the word obtained from $\bbu$ by deleting $u_{i}$, if $t=del$, or inserting an $x$ after $u_i$, if $t=ins_x$\footnote{Here, ``inserting after $u_0$'' means inserting before $u_1$.}. Moreover, we denote $\tilde{t}$ as the insertion/deletion types of $t$, where 
$$\tilde{t}=\begin{cases}
    del,~\text{if}~t=del;\\
    ins,~\text{if}~t=ins_x~\text{for some $x\in [q]$.}
\end{cases}$$

For non-negative integers $n$ and $l$, let $I=(i_l,i_{l-1},\ldots,i_1)$ be a non-increasing sequence of non-negative integers $n\geq i_l\geq i_{l-1}\geq \cdots \geq i_{1}\geq 0$, and let $T=(t_l,\ldots,t_1)\in \{del,ins_1,ins_2,\ldots,ins_q\}^{l}$ be a sequence of insertion/deletion operations, where $i_j\neq 0$ and $i_{j}>i_{j-1}$ if $t_j=del$. Then, we call the pair $(I,T)$ a sequence of $l$ insertions and deletions, and we write
\begin{equation}\label{eq_def_Phi}
    \Phi_{I,T}(\bbu)\eqdef(\phi_{i_1,t_1}\circ\phi_{i_2,t_2}\circ\cdots\circ\phi_{i_l,t_l})(\bbu)
\end{equation}
for the composition of the operations $\phi_{i_l,t_l}$ through $\phi_{i_1,t_1}$ applied to word $\bbu\in [q]^{n}$. Moreover, we denote $\tilde{T}=(\tilde{t}_l,\ldots,\tilde{t}_1)$ as the insertion/deletion type of $T$, and denote $(I,\tilde{T})$ as the insertion/deletion types of $(I,T)$.

As pointed out in \cite{ABGHK24}, whenever one obtains a word $\bbv$ from $\bbu$ by inserting and deleting several symbols, one can always reorder these operations to obtain sequence $(I,T)$ of insertion/deletion operations such that $\bbv=\Phi_{I,T}(\bbu)$. Since the elements of $I$ are non-increasing, an earlier operation cannot shift the location of a later operation. Thus, $i_j$ is not only the position in $(\phi_{i_{j+1},t_{j+1}}\circ\phi_{i_2,t_2}\circ\cdots\circ\phi_{i_l,t_l})(\bbu)$ at which the operation $\phi_{i_j,t_j}$ is applied, but also the original position in $\bbu$ where the operation occurs. This builds up a one-to-one correspondence between the word $\bbv$ obtained from $\bbu$ and the sequence $(I,T)$ satisfying $\bbv=\Phi_{I,T}(\bbu)$.

In \cite{ABGHK24}, Alon et al. proved the following lemma which bounds the number of $\lambda$-isolated elements in $I$ when $\bbu$ and $\bbv$ are $\lambda$-nonrepeating words with small edit distance.

\begin{lemma}\label{lem_bound_on_isolated_elements}(\cite{ABGHK24}, Lemma 6)
Let $n,d,\lambda\geq 1$, and let $\bbu,\bbv\in [q]^n$ be $\lambda$-nonrepeating words such that $\bbv=\Phi_{I,T}(\bbu)$ for some sequence of operations $(I,T)$. If the number of $\lambda$-isolated elements of $I$ is at least $d-1$, then $d_{E}(\bbu,\bbv)\geq d$.
\end{lemma}

Based on this result, we obtain the following bound on the number of edges in the neighbourhood of every vertex of graph $\Gamma_{n-\frac{d}{2}+1,\lambda}$ for $\lambda=3$. 

\begin{lemma}\label{lem_triangles_in_Gamma_alpha}
    Let $n$, $d$ be positive integers such that $6\leq d\leq \frac{n}{9}+2$ and let $\lambda=3$. 
    Let $q\geq n$ be a prime power and $\cC=RS_{\bm{\alpha}}(n,n-d/2+1)$ be an RS code over $\Fq$. Then, for every $\lambda$-nonrepeating word $\bbf\in\cC$, the number of pairs $(\bbg,\bbh)\in \cC^{2}$ such that both $\bbg,\bbh$ are $\lambda$-nonrepeating and $d_E(\bbf,\bbg)\leq d-2$, $d_E(\bbf,\bbh)\leq d-2$, $d_E(\bbg,\bbh)\leq d-2$ is at most
    $$2^{5(d-2)}\cdot d\cdot \left(\frac{\mathrm{e}n}{d-2}\right)^{\frac{3d-6}{2}}.$$
\end{lemma}


\begin{IEEEproof}
Let $k=n-d/2+1$. We call a triple $(\bbf,\bbg,\bbh)\in \cC^{3}$ good if $LCS(\bbf,\bbg)\geq k$, $LCS(\bbf,\bbh)\geq k$ and $LCS(\bbg,\bbh)\geq k$. Clearly, by (\ref{eq_ed_distance}), $d_E(\bbf,\bbg)\leq d-2$, $d_E(\bbf,\bbh)\leq d-2$ and $d_E(\bbg,\bbh)\leq d-2$ if and only if $(\bbf,\bbg,\bbh)$ is a good triple. 
Since $\cC$ is an $[n,k]$-RS code, each codeword $\bbf\in \cC$ is uniquely determined by $I$ and $\bbf|_{I}$ for any $I\subseteq [n]$ of size $k$. 
In the following, we use this property to derive an upper bound on the number of $\lambda$-nonrepeating good triples $(\bbf,\bbg,\bbh)\in \cC^{3}$ for a given $\lambda$-nonrepeating word $\bbf\in \cC$.




Note that a good triple $(\bbf,\bbg,\bbh)$ is uniquely determined by $\bbf$ and sequences $(I,T)$, $(I',T')$ of insertions and deletions for which $|I|=d_E(\bbf,\bbg)$, $\bbg=\Phi_{I,T}(\bbf)$ and $|I'|=d_E(\bbg,\bbh)$, $\bbh=\Phi_{I',T'}(\bbg)$. Since $\Phi_{I',T'}(\Phi_{I,T}(\bbf))=\bbh$, we can combine $(I,T)$ and $(I',T')$ to obtain a sequence $(I'',T'')$ of insertions and deletions of length $|I''|=|I|+|I'|\leq 2d-4=4(n-k)$ such that $\Phi_{I'',T''}(\bbf)=\bbh$. 

On the other hand, note that the insertion/deletion type $(I,\tilde{T})$ of $(I,T)$ actually specifies an LCS of $\bbf$ and $\bbg$ and the positions where it lies in both $\bbf$ and $\bbg$. Similarly, $(I',\tilde{T}')$ specifies an LCS of $\bbg$ and $\bbh$ and the positions where it lies in both $\bbg$ and $\bbh$. Since $LCS(\bbf,\bbg)\geq k$ and $LCS(\bbg,\bbh)\geq k$, once $\bbf$ is fixed in a good triple $(\bbf,\bbg,\bbh)$, $\bbg$ and $\bbh$ are uniquely determined by $(I,\tilde{T})$ and $(I',\tilde{T'})$, respectively. Next, we estimate the number of insertion/deletion types $(I'',\tilde{T}'')$ of all the insertions and deletions sequences $(I'',T'')$ satisfying $|I''|\leq 2d-4$, $\bbh=\Phi_{I'',T''}(\bbf)$ is $\lambda$-nonrepeating and $LCS(\bbf,\bbh)\geq k$. Note that  for a given $\bbf$, the number of good triples $(\bbf,\bbg,\bbh)$ is at most $2^{2d-4}$ times this number, as there are at most $2^{2d-4}$ choices of $(I,\tilde{T})$ and $(I',\tilde{T}')$ that produce the same sequence $(I'',\tilde{T}'')$. 

By Lemma \ref{lem_bound_on_isolated_elements}, $LCS(\bbf,\bbh)\geq k$, i.e., $d_E(\bbf,\bbh)\leq d-2$ implies that at most $d-2=2(n-k)$ elements of $I''$ are $\lambda$-isolated. Next, we bound the total number of ways to pick such an $I''$. Indeed, we define an equivalence relation $\sim$ on the elements of $I''$ by setting $i\sim j$ if $|i-j|\leq 2\lambda$ and then taking the transitive closure. Let $Q$ denote the number of equivalence classes. There are $Q_1\leq d-2$ equivalence classes of size $1$ coming from $\lambda$-isolated elements. There are $Q_2\leq 2$ equivalence classes coming from boundary elements $i\in I''$ satisfying $i\leq \lambda$ or $i\geq n-\lambda$, which we call \emph{boundary} equivalence classes. Note that there are at most $|I''|-(Q_1+Q_2)$ non-boundary and non-$\lambda$-isolated elements in $I''$ and each of these elements is not $\lambda$-isolated. Hence, by $|I''|\leq 2d-4$, we have
\begin{align*}
    Q&\leq \frac{|I''|-(Q_1+Q_2)}{2}+(Q_1+Q_2)=\frac{|I''|+Q_1+Q_2}{2}\\
    &\leq \frac{3}{2}(d-2)+1
\end{align*}
and $Q-Q_2\leq \frac{|I''|+Q_1}{2}\leq \frac{3}{2}(d-2)$. Moreover, there are at most ${n\choose Q-Q_2}$ ways to choose the minimal elements of the non-boundary equivalence classes, $\lambda^{Q_2}$ ways to choose the minimal elements of the boundary equivalence classes, and then $(2\lambda)^{|I''|-Q}$ ways to choose the remaining elements of $I''$.
Note that 
\begin{align*}
    & {n\choose Q-Q_2}\cdot \lambda^{Q_2}\cdot (2\lambda)^{|I''|-Q} \\
    \overset{(i)}{\leq} & (2\lambda)^{|I''|-(Q-Q_2)}\cdot \left(\frac{\mathrm{e}n}{Q-Q_2}\right)^{Q-Q_2}\\
    \overset{(ii)}{= } &  6^{|I''|}\cdot \left(\frac{\mathrm{e}n}{6(Q-Q_2)}\right)^{Q-Q_2}\\
    \overset{(iii)}{\leq} &  6^{2d-4}\cdot \left(\frac{\mathrm{e}n}{9(d-2)}\right)^{\frac{3d-6}{2}}\\
    = &\left(\frac{2}{\sqrt{3}}\right)^{2d-4}\cdot \left(\frac{\mathrm{e}n}{d-2}\right)^{\frac{3d-6}{2}},
\end{align*}
where (i) follows  by Stirling's formula, ${n\choose Q-Q_2}\leq \left(\frac{\mathrm{e}n}{Q-Q_2}\right)^{Q-Q_2}$; (ii) follows from $2\lambda=6$; (iii) follows from $|I''|\leq 2d-4$, $Q-Q_2\leq \frac{3(d-2)}{2}\leq \frac{n}{6}$, and $\left(\frac{\mathrm{e}n}{6 x}\right)^{x}$ is an increasing function of $x$ when $0<x\leq \frac{n}{6}$.\footnote{One can easily check that the derivative of $\left(\frac{\mathrm{e}n}{6 x}\right)^{x}$ equals $\left(\frac{\mathrm{e}n}{6 x}\right)^{x}(\ln(\frac{\mathrm{e}n}{6 x})-1)$, which implies $\left(\frac{\mathrm{e}n}{6 x}\right)^{x}$ is increasing in $x$ when $0<x\leq \frac{n}{6}$.}
Since there are at most $(\frac{3}{2}(d-2)+1)\cdot 3 \leq 5d$ possible choices of $Q$ and $Q_2$,  the total number of ways to choose $I''$ is at most $5d\cdot \left(\frac{2}{\sqrt{3}}\right)^{2d-4}\cdot \left(\frac{\mathrm{e}n}{d-2}\right)^{\frac{3d-6}{2}}$. 

Since $\tilde{T}''\in \{del,ins\}^{|I''|}$ and $\abs{I''}\leq 2d-4$, there are at most $2^{2d-4}$ ways to pick $\tilde{T}''$. Therefore, the total number of insertion/deletion types $(I'',\tilde{T}'')$ of all the insertions and deletions sequences $(I'',T'')$ satisfying $|I''|\leq 2d-4$, $\bbh=\Phi_{I'',T''}(\bbf)$ is $\lambda$-nonrepeating and $LCS(\bbf,\bbh)\geq k$ is at most 
\begin{align*}
2^{2d-4}\cdot 5d\cdot \left(\frac{2}{\sqrt{3}}\right)^{2d-4}\cdot \left(\frac{\mathrm{e}n}{d-2}\right)^{\frac{3d-6}{2}}\leq 2^{3(d-2)}\cdot d\cdot \left(\frac{\mathrm{e}n}{d-2}\right)^{\frac{3d-6}{2}}.
\end{align*}
Thus, the number of $\lambda$-nonrepeating good triples $(\bbf,\bbg,\bbh)$ for a given $\lambda$-nonrepeating word $\bbf\in \cC$ is at most $2^{2(d-2)}\cdot 2^{3(d-2)}\cdot d\cdot \left(\frac{\mathrm{e}n}{d-2}\right)^{\frac{3d-6}{2}}$, which concludes the proof.
\end{IEEEproof}


The proof of (\ref{eq_GV_type_lower_bound}) in Theorem \ref{Thm_GV_type_lower_bound} is now immediate.

\begin{IEEEproof}[Proof of (\ref{eq_GV_type_lower_bound}) in Theorem \ref{Thm_GV_type_lower_bound}]
   Consider the graph $\Gamma_{n-\frac{d}{2}+1,\lambda}$ with $\lambda=3$. By Lemma \ref{lem_bound_on_lambda-repaeating_words}, we have $|V(\Gamma_{n-\frac{d}{2}+1,\lambda})|\geq q^{n-\frac{d}{2}+1}(1-\frac{1}{n})$. By Lemma \ref{lem_triangles_in_Gamma_alpha}, for every $\bbu\in V(\Gamma_{n-\frac{d}{2}+1,\lambda})$, it's neighborhood contains at most ${2^{5(d-2)}\cdot d\cdot \left(\frac{\mathrm{e}n}{d-2}\right)^{\frac{3d-6}{2}}}$ edges. Moreover, by (\ref{eq_maximum_degree_Gamma_{k,lambda}}), we have $\Delta\parenv{\Gamma_{n-\frac{d}{2}+1,\lambda}}< {n\choose \frac{d}{2}-1}^2$.  Set $\Delta={n\choose \frac{d}{2}-1}^2$ and $f=\Delta^2/\left({2^{5(d-2)}\cdot d\cdot \left(\frac{\mathrm{e}n}{d-2}\right)^{\frac{3d-6}{2}}}\right)$. 
   Then, the condition of Lemma \ref{lem_chi_bounded_neighborhood} is satisfied, and so, from (\ref{eq_lb_independent_bounded_neighborhood}) we have 
   $$D_q(n,d)\geq \alpha(\Gamma_{n-\frac{d}{2}+1,\lambda})\geq   \ln(\min\{f,\Delta\}) \cdot \frac{q^{n-\frac{d}{2}+1}}{{n\choose \frac{d}{2}-1}^2}\cdot(1-o(1)).$$

  By ${n\choose \frac{d}{2}-1}\geq \left(\frac{n}{\frac{d}{2}-1}\right)^{\frac{d-2}{2}}$, we have $\ln \Delta \geq \frac{d-2}{2}(\ln n - \ln (d-2)+\ln 2)$ and $\Delta^2\geq 2^{2d-4}\cdot \left(\frac{n}{d-2}\right)^{2d-4}$.
   The latter implies that
   \begin{align*}
       f &\geq \left(\frac{n}{2^6\mathrm{e}^3(d-2)}\right)^{\frac{d-2}{2}}\cdot \frac{1}{d}. 
   \end{align*}
   Since $\ln(2^6\mathrm{e}^3)<7.2$ and $\frac{2\ln{d}}{d-2}<0.9$ for $d\geq 6$, we have 
   \begin{align*}
       \ln f &= (\frac{d}{2}-1)(\ln n - \ln (d-2)-\ln(2^6\mathrm{e}^3))-\ln{d}\\
       &\geq (\frac{d}{2}-1)(\ln n - \ln (d-2)-8.1)
   \end{align*}
   Hence, $\ln(\min\{f,\Delta\})\geq (\frac{d}{2}-1)(\ln{n}-\ln{(d-2)}-8.1)$. This concludes the proof.
\end{IEEEproof}

When $k=n$, by definition, we have $\Gamma_{n}=\Gamma$. Thus, for $1\leq \lambda\leq n$, the graph $\Gamma_{n,\lambda}$ is actually the subgraph of $\Gamma$ induced by all the $\lambda$-nonrepeating words in $[q]^n$. Note that for $\lambda\geq \lceil3\log_q{n}\rceil$, the total number of $\lambda$-repeating words in $[q]^n$ is at most ${n\choose 2}q^{n-\lambda}=o_n(q^n)$. Thus, similarly, we can apply Lemma \ref{lem_chi_bounded_neighborhood} and (\ref{eq_lb_independent_bounded_neighborhood}) to $\Gamma_{n,\lambda}$ and obtain the lower bound (\ref{eq_GV_type_lower_bound2}) in Theorem \ref{Thm_GV_type_lower_bound}. Since the proof of (\ref{eq_GV_type_lower_bound2}) uses Lemma \ref{triangles_in_Gamma_Alon} and is similar to that of (\ref{eq_GV_type_lower_bound}). Thus, we left it in the appendix for interested readers.








\appendices

\section{Proof of (\ref{eq_GV_type_lower_bound2}) in Theorem \ref{Thm_GV_type_lower_bound}}\label{sec_appendix_a}

For completeness, we first include the following proof of Lemma \ref{triangles_in_Gamma_Alon}, where we demonstrate all the exact forms of the constant factors that were omitted in the proof of Lemma \ref{triangles_in_Gamma_Alon_orin} in \cite{ABGHK24}.

\begin{IEEEproof}[Proof of Lemma \ref{triangles_in_Gamma_Alon}]
We call a triple $(\bbu,\bbv,\bbw)\in([q]^n)^3$ \emph{good} if $d_E(\bbu,\bbv)\leq a$, $d_E(\bbv,\bbw)\leq b$ and $d_E(\bbu,\bbw)\leq c$. Note that $d_E(\bbu,\bbv)\leq d_E(\bbv,\bbw)+d_E(\bbw,\bbu)\leq b+c$ by the triangle inequality, so all good triples $(\bbu,\bbv,\bbw)$ satisfy $d(\bbu,\bbv)\leq b+c$. Thus, we restrict our attention to the regime $a,b,c\geq 2$ and $a\leq b+c$.

Observe that the probability of a uniformly random $\bbu\in[q]^n$ being $\lambda$-repeating is at most ${n\choose 2}q^{-\lambda}$. Thus, the total number of $\lambda$-repeating words in $[q]^n$ is at most ${n\choose 2}q^{n-\lambda}$. By Lemma \ref{bounds_on_deletion/insertion_ball}, for each $\bbu\in [q]^n$, there are $n^{2a}q^a(1+o_{n}(1))$
words $\bbv$ at distance at most $a$ and  $n^{2c}q^c(1+o_{n}(1))$
words $\bbw$ at distance at most $c$. Thus, there are at most $n^{2a+2c}q^{a+c}(1+o_{n}(1))$ good triples $(\bbu,\bbv,\bbw)$ for each choice of $\bbu$. Therefore, the total number of good triples containing a $\lambda$-repeating word is at most 
\begin{align*}
    3\cdot {n\choose 2}q^{n-\lambda}\cdot n^{2a+2c}q^{a+c}(1+o_{n}(1)). 
\end{align*}
When $\lambda=\lceil 10a\log_q{n}\rceil$, this implies that there are at most 
\begin{align}
    q^{n+2a}n^{-4a}(3+o_{n}(1)) \label{eq0_triangles_in_Gamma_Alon}
\end{align}
good triples containing a $\lambda$-repeating word. 

Next, we bound the number of good triples consisting of $\lambda$-nonrepeating word.

Note that a good triple $(\bbu,\bbv,\bbw)$ is uniquely determined by $\bbu$ and sequences $(I,T)$, $(I',T')$ of insertions and deletions for which $\bbw=\Phi_{I,T}(\bbu)$ and $\bbv=\Phi_{I',T'}(\bbw)$. Since $d_E(\bbu,\bbw)\leq c$ and $d_{E}(\bbw,\bbv)\leq b$, we may choose $(I,T)$ to have length at most $2c$ and $(I',T')$ to have length at most $2b$. Since $\Phi_{I',T'}(\Phi_{I,T}(\bbu))=\bbv$, we can ``combine'' the operations of $(I,T)$ and $(I',T')$ to obtain a sequence $(I'',T'')$ of insertions and deletions of length $|I''|=|I|+|I'|\leq 2b+2c$ such that $\Phi_{I'',T''}(\bbu)=\bbv$. Furthermore, there are at most ${|I''|\choose |I|}\leq 2^{2b+2c}$ choices of $(I,T)$ and $(I',T')$ that produce each such sequence $(I'',T'')$. Thus, for a given $\bbu$, the number of good triples $(\bbu,\bbv,\bbw)$ is at most $2^{2b+2c}$ times the number of ways to pick a sequence $(I'',T'')$ such that $|I''|\leq 2b+2c$, $\bbv=\Phi_{I'',T''}(\bbu)$ is $\lambda$-nonrepeating and $d_{E}(\bbu,\bbv)\leq a$.

By Lemma \ref{lem_bound_on_isolated_elements}, the assumption $d_{E}(\bbu,\bbv)\leq a$ implies that at most $2a$ of the elements of $I''$ are $\lambda$-isolated. We claim that the total number of ways to pick such an $I''$ is at most $n^{a+b+c}(22a)^{2b+2c+1}(\log{n})^{2b+2c}$ when $\lambda=\lceil 10a\log_q{n}\rceil$. Then, since there are $(q+1)^{2b+2c}\leq (\frac{3q}{2})^{2b+2c}$ ways to pick $T''$, the number of good triples consisting of $\lambda$-nonrepeating words for fixed $\bbu$ is at most
\begin{align}
    n^{a+b+c}q^{2b+2c}(66a)^{2b+2c+1}(\log{n})^{2b+2c}.\label{eq1_triangles_in_Gamma_Alon}
\end{align}

To prove the claim, we define an equivalence relation $\sim$ on elements of $I''$ by setting $i\sim j$ if $|i-j|\leq 2\lambda$ and then taking the transitive closure. Let $Q$ denote the number of equivalence classes. There are $Q_1\leq 2a$ equivalence classes of size $1$ coming from $\lambda$-isolated elements, and $Q_2\leq 2$ equivalence classes coming from elements $i\in I''$ satisfying $i\leq \lambda$ or $i\geq n-\lambda$, which we call boundary equivalence classes. Hence, we have
\begin{align*}
    Q&\leq (|I''|-(Q_1+Q_2)/2)+Q_1+Q_2=(|I''|+Q_1+Q_2)/2\\
    &\leq a+b+c+1
\end{align*}
and $Q-Q_2\leq (|I''|+Q_1)/2\leq a+b+c$. Moreover, there are at most ${n\choose Q-Q_2}$ ways to choose the minimal elements of the non-boundary equivalence classes, $\lambda^{Q_2}$ ways to choose the minimal elements of the boundary equivalence classes, and then $(2\lambda)^{|I''|-Q}$ ways to choose the remaining elements of $I''$. Thus, the total number of ways to choose $I''$ is at most ${n\choose Q-Q_2}\cdot \lambda^{Q_2}\cdot (2\lambda)^{|I''|-Q}$ multiplying by $3(a+b+c+1)$ (for the possible choices of $Q$ and $Q_2$), which is upper bound by
\begin{align}
    &3(a+b+c+1)\cdot {n\choose Q-Q_2}\cdot \lambda^{Q_2}\cdot (2\lambda)^{|I''|-Q} \nonumber\\
    \overset{(i)}{\leq} & 12a\cdot {n\choose Q-Q_2}\cdot (2\lambda)^{|I''|-(Q-Q_2)} \nonumber\\
    \overset{(ii)}{\leq} & 12a \cdot \parenv{\frac{\mathrm{e}n}{2\lambda(Q-Q_2)}}^{Q-Q_2}\cdot (2\lambda)^{|I''|}, \label{eq2_triangles_in_Gamma_Alon}
\end{align}
where (i) follows from $3(a+b+c+1)\leq 12a$ and $\lambda^{Q_2}\leq (2\lambda)^{Q_2}$; (ii) follows from ${n\choose Q-Q_2}\leq \parenv{\frac{\mathrm{e}n}{Q-Q_2}}^{Q-Q_2}$ by Stirling's formula. When $\lambda=\lceil 10a\log_q{n}\rceil$, by $\parenv{\frac{\mathrm{e}n}{2\lambda(Q-Q_2)}}^{Q-Q_2}\leq n^{Q-Q_2}$ and $Q-Q_2\leq a+b+c$, (\ref{eq1_triangles_in_Gamma_Alon}) is upper bounded by 
\begin{align*}
    12a \cdot n^{a+b+c}\cdot (2\lambda)^{|I''|}
    \leq & n^{a+b+c}(22a)^{2b+2c+1}(\log{n})^{2b+2c},
\end{align*}
where the inequality follows from $\lambda\leq 11a\log{n}$ and $12a\leq 22a$. This confirms the claim.

Finally, recall that the number of $\lceil 10a\log_q{n}\rceil$-repeating words is at most $q^nn^{-8a}=o_{n}(q^n)$, thus there are $q^n(1-o_{n}(1))$ different choices of $\lceil 10a\log_q{n}\rceil$-nonrepeating word $\bbu$. Then, by (\ref{eq0_triangles_in_Gamma_Alon}) and (\ref{eq1_triangles_in_Gamma_Alon}), we can conclude the proof.
\end{IEEEproof}

To use Lemma \ref{lem_chi_bounded_neighborhood} and (\ref{eq_lb_independent_bounded_neighborhood}), we need to upper bound the number of edges within the neighborhood of any vertex $\bbu\in V(\Gamma_{n,\lambda})$. By taking $a=\frac{d-1}{2}$ in the proof of (\ref{eq1_triangles_in_Gamma_Alon}) in Lemma \ref{triangles_in_Gamma_Alon}, we have the following result.

\begin{corollary}\label{coro_neighbothood_Gamma}
    Let $\lambda=\lceil5(d-2)\log_q{n}\rceil$. Then, for every $\bbu \in V(\Gamma_{n,\lambda})$, we have
    $$N_{\Gamma_{n,\lambda}}(\bbu)\leq c_{d,q}\cdot n^\frac{3(d-2)}{2}(\log{n})^{2(d-2)},$$
    where $c_{d,q}$ is a constant depends on $q$ and $d$.
\end{corollary}

The proof of (\ref{eq_GV_type_lower_bound2}) in Theorem \ref{Thm_GV_type_lower_bound} is now immediate.

\begin{IEEEproof}[Proof of (\ref{eq_GV_type_lower_bound2}) in Theorem \ref{Thm_GV_type_lower_bound}]
    Take $\lambda=\lceil5(d-2)\log_q{n}\rceil$, we have $|V(\Gamma_{n,\lambda})|= q^{n}(1-{n\choose 2}q^{-\lambda})=q^n(1-o_n(1))$. By Corollary \ref{coro_neighbothood_Gamma}, for every $\bbu\in V(\Gamma_{n,\lambda})$, the number of edges in it's neighborhood is at most $c_{d,q}\cdot n^\frac{3(d-2)}{2}(\log{n})^{2(d-2)}$. Moreover, by Lemma \ref{bounds_on_deletion/insertion_ball}, we have $\Delta\parenv{\Gamma_{n,\lambda}}\leq \Delta\parenv{\Gamma}\leq n^{d-2}q^{\frac{d}{2}-1}(1+o_n(1))$. Then, applying Lemma \ref{lem_chi_bounded_neighborhood} and (\ref{eq_lb_independent_bounded_neighborhood}) on $\Gamma_{n,\lambda}$ with $f=\frac{n^{\frac{d}{2}-1}}{c_{d,q}\cdot q^{d-2}(\log{n})^{2d-4}}$, we have $$D_q(n,d)\geq \alpha(\Gamma_{n,\lambda})\geq q^{n-\frac{d}{2}+1}\cdot \frac{\log{n}}{n^{d-2}}\cdot (\frac{d}{2}-1+o_{n}(1)).$$
    This concludes the proof.
\end{IEEEproof}


\bibliographystyle{IEEEtran}
\bibliography{biblio}

\end{document}